%% file: submission_to_AnnProbab.tex
\documentclass[a4paper,reqno]{amsart}
\usepackage[foot]{amsaddr}
\usepackage{fullpage}
\usepackage{graphicx}
\usepackage{color}
\usepackage{amsthm,amssymb}
\usepackage{xfrac}
\usepackage{subcaption}

\usepackage{tikz}
\usepackage{tikz-cd}
\usepackage[hidelinks]{hyperref}
\usepackage{mathtools}
\usepackage{xcolor}
\usepackage{algpseudocode}
\usepackage{algorithm}
\usepackage{enumerate}
\usepackage{nicematrix}
\usepackage{adjustbox}
\usepackage[normalem]{ulem}

\usepackage{lmodern}

\usepackage[colorinlistoftodos,prependcaption,textsize=tiny]{todonotes}

\newcommand {\mm}[1]   {\ifmmode{#1}\else{\mbox{\(#1\)}}\fi}

\newcommand{\Expect}[1]     {\mm{{\mathbb E}\left[{#1}\right]}}
\newcommand{\Prob}[1]       {\mm{{\mathbb P}[{#1}]}}
\newcommand{\Length}[1]     {\mm{\|{#1}\|}}
\newcommand{\EMST}[1]       {\mm{{\rm EMST}{({#1})}}}
\newcommand{\Cmst}          {\mm{c}}
\newcommand{\CLmst}         {\mm{c_L}}

\newcommand{\Rspace}        {\mm{{\mathbb R}}}

\newcommand{\Hgroup}[1]     {\mm{\sf H}_{#1}}

\newcommand{\kkk}           {\mm{{\bf k}}}

\newcommand{\Domain}[2]     {\mm{\rm Dom}_{#1}}
\newcommand{\Image}[2]      {\mm{\rm Img}_{#1}}
\newcommand{\Codomain}[2]   {\mm{\rm Cod}_{#1}}
\newcommand{\Cokernel}[2]   {\mm{\rm Cok}_{#1}}
\newcommand{\Relative}[2]   {\mm{\rm Rel}_{#1}}
\newcommand{\Kernel}[2]     {\mm{\rm Ker}_{#1}}

\newcommand{\Ckernel}[1]    {\mm{C^{\tt ker}_{#1}}}
\newcommand{\Cimage}[1]     {\mm{C^{\tt img}_{#1}}}
\newcommand{\Ccokernel}[1]  {\mm{C^{\tt cok}_{#1}}}
\newcommand{\Cdomain}[1]    {\mm{C^{\tt dom}_{#1}}}
\newcommand{\Ccodomain}[1]  {\mm{C^{\tt cod}_{#1}}}
\newcommand{\Crelative}[1]  {\mm{C^{\tt rel}_{#1}}}

\newcommand{\Delaunay}[1]   {\mm{{\rm Del}{({#1})}}}

\newcommand{\Radiusf}       {\mm{{f}}}

\newcommand{\rank}[1]       {\mm{{\rm rank\,}{#1}}}

\newcommand{\diff}          {\mm{\,}{\rm d}}

\newcommand{\Dshort}        {\mm{D^{\rm sh}}}
\newcommand{\Dlong}         {\mm{D^{\rm lg}}}
\newcommand{\Bcertified}    {\mm{B^{\rm ct}}}
\newcommand{\Buncertified}  {\mm{B^{\rm un}}}

\newcommand{\norm}[1]       {\mm{\|{#1}\|}}

\newcommand{\Skip}[1]       {}

\definecolor{blue-red}{rgb}{0.8, 0.00, 0.95}

\theoremstyle{definition}

\newtheorem{theorem}{Theorem}
\numberwithin{theorem}{section}

\newtheorem{lemma}[theorem]{Lemma}

\numberwithin{equation}{section}

\title{Expected Length of the Euclidean Minimum Spanning Tree and 1-norms of Chromatic Persistence Diagrams in the Plane}

\author{Ond\v{r}ej Draganov$^{1,5}$}
\email{$^5$ondrej.draganov@inria.fr}
\author{Herbert Edelsbrunner$^{2,6}$}
\email{$^6$edels@ist.ac.at}
\author{Sophie Rosenmeier$^{3,4,7}$}
\email{$^7$rosenmeier@biochem.mpg.de}
\author{Morteza Saghafian$^{2,8}$}
\email{$^8$morteza.saghafian@ist.ac.at}

\address{$^{1}$INRIA, Sophia Antipolis, France}
\address{$^{2}$ISTA (Institute of Science and Technology Austria), Kloster\-neu\-burg, Austria}
\address{$^{3}$University of Vienna, Vienna, Austria}
\address{$^{4}$Max Planck Institute of Biochemistry, Munich, Germany}

\thanks{\emph{Funding.} Subsets of the authors are partially supported by the European Research Council (ERC) Horizon 2020 project `Alpha Shape Theory Extended', grant no.\ 788183, partially by the DFG Collaborative Research Center TRR 109, `Discretization in Geometry and Dynamics', Austrian Science Fund (FWF), grant no.\ I 02979-N35, and by the Agence Nationale de la Recherche project `AI4scMed', France 2030 ANR-22-PESN-000.}


\begin{document}

\begin{abstract}
  Let $\Cmst$ be the constant such that the expected length of the Euclidean minimum spanning tree of $n$ random points in the unit square is $\Cmst \sqrt{n}$ in the limit, when $n$ goes to infinity.
  We improve the prior best lower bound of $0.6008 \leq \Cmst$ by Avram and Bertsimas~\cite{AvBe92} to $0.6289 \leq \Cmst$.
  The proof is a by-product of studying the persistent homology of randomly $2$-colored point sets.
  Specifically, we consider the filtration induced by the inclusions of the two mono-chromatic sublevel sets of the Euclidean distance function into the bi-chromatic sublevel set of that function.
  Assigning colors randomly, and with equal probability, we show that the expected $1$-norm of each chromatic persistence diagram is a constant times $\sqrt{n}$ in the limit,
  and we determine the constant in terms of $\Cmst$ and another constant, $\CLmst$, which arises for a novel type of Euclidean minimum spanning tree of $2$-colored point sets.
\end{abstract}

\maketitle

\section{Introduction}
\label{sec:1}

A classic question in stochastic geometry is the expected length of the minimum spanning tree of $n$ points sampled uniformly at random in the unit square.
As proved by Steele~\cite{Ste88}, there exists a constant $\Cmst > 0$ such that the expected length is $\Cmst \sqrt{n}$ in the limit, when $n$ goes to infinity; but see also Beardwood, Halton and Hammersley~\cite{BHH59}.
The constant is however not known.
The current best upper bound is $\Cmst \leq 0.7072$ due to Gilbert~\cite{Gil65}, and the best lower bound prior to this paper is $0.6008 \leq \Cmst$ due to Avram and Bertsimas~\cite{AvBe92}; see also \cite[page 50]{Yuk98}.
We improve the lower bound to $0.6289 \leq \Cmst$.

\smallskip
Our interest in the length of the Euclidean minimum spanning tree, or EMST for short, originates in its connection to the $1$-norms of the chromatic persistence diagrams of the Euclidean distance function, which were recently introduced as a tool in topological data analysis \cite{CDES23}.
Specifically, \cite{CDES24} defines the \emph{MST-ratio} of a pair of point sets, $B \subseteq A$, as the total length of the EMSTs of $B$ and $A \setminus B$ divided by the length of the EMST of $A$, and proves tight upper and lower bounds for lattices in $\Rspace^2$.
Additional bounds including for points sampled uniformly at random from the unit square and in general metric spaces can be found in \cite{DPT23} and in \cite{AMS24}, respectively.
This paper goes beyond the MST-ratio and studies the expected $1$-norms of all six chromatic persistence diagrams.
More precisely, we consider the $2$-colored case, the inclusion maps from the sublevel sets of the two mono-chromatic distance functions to the sublevel set of the bi-chromatic distance function, and the six diagrams that measure the persistent homology of the domain, codomain, pair, image, kernel, and cokernel.
The $1$-norm of the degree-$0$ codomain diagram is half the length of the EMST of the points, and we establish a similar if less direct connection to the degree-$1$ codomain diagram, which measures the homology of loops.
Using analytic results for Poisson point processes in \cite{ENR17}, we combine the two connections and improve the lower bound for $\Cmst$.

\smallskip
Our approach is predominantly geometric and revolves around the radius function on the \emph{Delaunay mosaic} (aka the Delaunay triangulation) of the points.
To define it, call a circle \emph{empty} if all points of a given set either lie on or outside the circle.
Given a locally finite set $A \subseteq \Rspace^2$, this mosaic is the polyhedral complex whose cells are the convex hulls of the subsets $Q \subseteq A$ for which there is an empty circle that passes through all points in $Q$ and no other points in $A$.
The mosaic consists of vertices, edges, and convex polygons, and if the maximum number of points on a common circle is three, then all polygons are triangles.
Writing $\Delaunay{A}$ for this mosaic, the \emph{radius function}, denoted $f \colon \Delaunay{A} \to \Rspace$, maps $Q$ to the radius of the smallest empty circle that passes through all points in $Q$.
For each $r \geq 0$, the collection of cells with radius at most $r$ is traditionally referred to as the \emph{alpha complex} of $A$ and $r$ \cite{EKS83}.
Importantly, it is homotopy equivalent to the corresponding sublevel set of the Euclidean distance function, which is the union of closed disks of radius $r$ centered at the points in $A$.

\smallskip
Of special significance are the \emph{critical cells} $Q \in \Delaunay{A}$ of $f$, for which the smallest circle that encloses all points of $Q$ is also the smallest empty such circle.
Assuming $A$ is generic\footnote{We call $A \subseteq \Rspace^2$ \emph{generic} if no four points lie on a common circle and no three lie on a circle whose center lies on a common line with two of the three points.}, a triangle is critical iff all three of its angles are less than $\frac{\pi}{2}$ (the triangle is acute), and an edge is critical iff the opposite angle inside each triangle that shares the edge is less than $\frac{\pi}{2}$ (the circle centered at its midpoint with radius half the length of the edge is empty).
The special role played by the critical cells becomes apparent when we increase $r$ continuously and monitor the corresponding filtration of alpha complexes and unions of disks.
This process induces an ordering on the simplices, and we construct the Delaunay mosaic incrementally by adding the simplices according to this order.
Whenever we add a critical cell, the homotopy type of the alpha complex changes, and because of the homotopy equivalence, so does the homotopy type of the union of disks.
For a critical edge, the homotopy type may change because the edge connects two connected components or it forms a new loop in one such component.
In the former case, we say the edge \emph{gives death} to one component, while in the latter case it \emph{gives birth} to a loop.
\begin{figure}[hbt]
  \centering
  \includegraphics[width=0.7\linewidth]{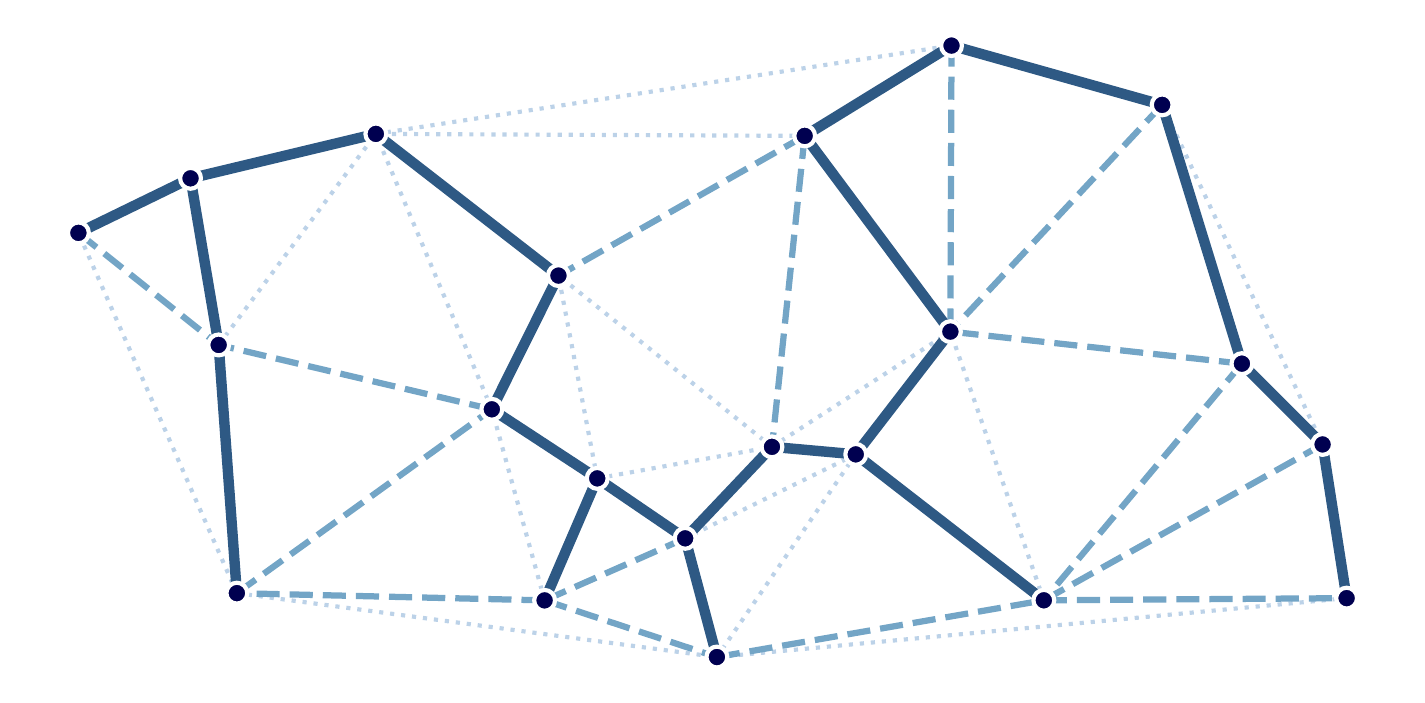}
  \vspace{-0.1in}
  \caption{\footnotesize A finite set of points in $\Rspace^2$, the Delaunay mosaic (\emph{all} edges), the Gabriel graph (\emph{dashed} and \emph{solid} edges), and the Euclidean minimum spanning tree (\emph{solid} edges).
  Except for the \emph{dotted} ones, \emph{all} edges are critical, the \emph{solid} ones give death, and the \emph{dashed} ones give birth.}
  \label{fig:delaunay_gabriel_mst}
\end{figure}
This categorization of the edges in the Delaunay mosaic can also be seen as the \emph{Gabriel graph} \cite{GaSo69}, which consists of all critical edges of the radius function, and the \emph{Euclidean minimum spanning tree} as its subset
of the death-giving edges; see Figure~\ref{fig:delaunay_gabriel_mst}.
As we will see shortly, the relative frequency of the different types of edges and triangles is at the core of all questions studied in this paper.
For a stationary Poisson point process in $\Rspace^2$, these frequencies become densities, which have been considered more than half a century ago by Miles~\cite{Mil70}.
For example, he showed that two thirds of the edges and half of the triangles in the Delaunay mosaic are critical in expectation.
However, in this paper we need the stronger results proved in \cite{ENR17}, which give the fractions for each $r \geq 0$.

\smallskip
The outline of this paper is as follows.
Section~\ref{sec:2} presents background on Delaunay mosaics for Poisson point processes in the plane.
Section~\ref{sec:3} focuses on the EMST and proves the new lower bound for the asymptotic constant of its expected length.
Section~\ref{sec:4} introduces the framework of chromatic persistence and mentions the generalization of the EMST whose asymptotic constant is $\CLmst$.
Section~\ref{sec:5} describes what we know about the expected $1$-norms of the six chromatic persistence diagrams in the $2$-colored case.
Section~\ref{sec:6} complements the analytic results with computational experiments and gives estimates for the two asymptotic constants, $\Cmst$ and $\CLmst$.
Section~\ref{sec:7} concludes the paper.

\section{Poisson--Delaunay Mosaics}
\label{sec:2}

In this section, we collect results about the Delaunay mosaic of a stationary Poisson point process in the plane, in particular about the expected moments of its critical edges and triangles, and the subsets of critical edges and triangles whose smallest enclosing circles cross the boundary of the unit square.

\subsection{Critical Edges and Triangles}
\label{sec:2.1}

To prepare for the improvement of the lower bound for the asymptotic constant, $c$, we introduce the expected moments of the radii of the critical edges and triangles in the Delaunay mosaic of a stationary Poisson point process in $\Rspace^2$.
We note that the Delaunay mosaic of such a process is generic with probability $1$.
In this case, an edge or triangle is critical iff the smallest empty circle that passes through its vertices is also the smallest enclosing circle.
For a threshold $r_0 \geq 0$, we write $N_1(r_0)$, $F_1(r_0)$, and $S_1(r_0)$ for the sums of the $0$-th, $1$-st, and $2$-nd powers\footnote{So $N_1(r_0)$ is the number of such critical edges, $F_1(r_0)$ is the sum of their radii (half-lengths), and $S_1(r_0)$ is the sum of their squared radii.} of the radii of all critical edges whose midpoints lie in $[0,1]^2$ and whose radii are less than or equal to $r_0$.
Similarly, we write $N_2(r_0)$, $F_2(r_0)$, and $S_2(r_0)$ for the $0$-th, $1$-st, and $2$-nd powers of the radii of all critical triangles whose circumcenters lie in $[0,1]^2$ and whose radii are less than or equal to $r_0$.
\begin{lemma}
  \label{lem:moments}
  Let $A$ be a stationary Poisson point process with intensity $n > 0$ in $\Rspace^2$, and $\Delaunay{A}$ its Delaunay mosaic.
  Then for every $r_0 > 0$, we have
  \begin{align}
    \Expect{N_1(r_0)} &= 2n \cdot \gamma(1, x);
      \hspace{0.25in}
    \Expect{F_1(r_0)} = \tfrac{2 \sqrt{n}}{\sqrt{\pi}} \cdot \gamma(\tfrac{3}{2}, x);
      \hspace{0.25in}
    \Expect{S_1(r_0)} = \tfrac{2}{\pi} \cdot \gamma(2, x) ; 
      \label{eqn:L21edges} \\
    \Expect{N_2(r_0)} &= n \cdot \gamma(2, x); 
      \hspace{0.32in}
    \Expect{F_2(r_0)}  = \tfrac{\sqrt{n}}{\sqrt{\pi}} \cdot \gamma(\tfrac{5}{2}, x);
      \hspace{0.30in}
    \Expect{S_2(r_0)}  = \tfrac{1}{\pi} \cdot \gamma(3, x) .
      \label{eqn:L21triangles}
  \end{align}
  in which $x = n \pi r_0^2$ and $\gamma(k,x) = \int_{t=0}^x t^{k-1} e^{-t} \diff t$ is the incomplete gamma function.
\end{lemma}
\begin{proof}
  We derive the first relations in \eqref{eqn:L21edges} and \eqref{eqn:L21triangles} from Theorem~1 in \cite{ENR17}, which is formulated for a stationary Poisson point process with intensity $\varrho > 0$ in $\Rspace^d$ and a Borel set $\Omega$ with $d$-dimensional volume $\norm{\Omega}$, and it gives expected numbers of intervals of simplices in $\Delaunay{A}$ of various types.
  An interval of type `$\ell \leq k$'
  in $\Delaunay{A}$ consists of a minimal simplex of dimension $\ell$, a maximal simplex of dimension $k$, and all simplices that contain the minimal and are contained in the maximal simplex, such that all these simplices share the same radius.
  For each type, the theorem counts the intervals that are maximal with respect to inclusion and for which the center of the smallest circumsphere of the maximal simplex lies in $\Omega$.
  The lemma asserts that the expected number of such intervals of type $\ell\leq k$, in $\Rspace^d$, whose radii are smaller than or equal to a threshold $r_0 > 0$ is
  \begin{align}
    \Expect{N_{\ell \leq k \leq d} (r_0)} &= \frac{C_{\ell \leq k \leq d} \cdot \varrho \norm{\Omega}}{\Gamma (k)} \cdot  \gamma (k, \varrho \nu_d r_0^d),
      \label{eqn:ENRN}
  \end{align}
  in which $\nu_d$ is the volume of the unit ball in $\Rspace^d$, $\Gamma(k) = (k-1)!$, and the $C_{\ell \leq k \leq d}$ are constants.
  To count the critical edges and critical triangles, we set $\ell = k = 1$ and $\ell = k = 2$, respectively.
  Further setting $d=2$ with $\nu_2 = \pi$, $\varrho = n$, and $\Omega = [0,1]^2$ with $\norm{\Omega} =1$, we get the first relations in \eqref{eqn:L21edges} and \eqref{eqn:L21triangles} because $\Gamma(1) = \Gamma(2) = 1$ and $C_{1 \leq 1 \leq 2} = 2$, $C_{2 \leq 2 \leq 2} = 1$, as proved in \cite{ENR17}.
  
  To get the remaining relations in \eqref{eqn:L21edges} and \eqref{eqn:L21triangles}, we add one and then two powers of the radius inside the integrals in the proof of \eqref{eqn:ENRN}.
  In particular, setting $r(t) = \sqrt[d]{t/(\varrho \nu_d)}$ and $x = \varrho \nu_d r_0^d$, the expected $1$-st and $2$-nd moments are
  \begin{align}
    \Expect{F_{\ell \leq k \leq d} (r_0)} &= \frac{C_{\ell \leq k \leq d}}{\Gamma(k)} \cdot \varrho \norm{\Omega} \cdot \int_{t=0}^x r(t) t^{k-1} e^{-t} \diff t 
    = \frac{C_{\ell \leq k \leq d} \cdot \varrho \norm{\Omega}}{\Gamma(k) \cdot \sqrt[d]{\varrho \nu_d}}
    \cdot \gamma(k+\tfrac{1}{d}, x) 
      \label{eqn:ENRF} \\
    \Expect{S_{\ell \leq k \leq d} (r_0)} &= \frac{C_{\ell \leq k \leq d}}{\Gamma(k)} \cdot \varrho \norm{\Omega} \cdot \int_{t=0}^x r^2(t) t^{k-1} e^{-t} \diff t 
    = \frac{C_{\ell \leq k \leq d} \cdot \varrho \norm{\Omega}}{\Gamma(k) \cdot \sqrt[d/2]{\varrho \nu_d}} \cdot
    \gamma(k+\tfrac{2}{d}, x) .
      \label{eqn:ENRS}
  \end{align}
  Setting again $d=2$ with $\nu_2 = \pi$, $\varrho = n$, and $\Omega = [0,1]^2$, as well as $\ell = k = 1$ and $\ell = k = 2$, and recalling $C_{1 \leq 1 \leq 2} = 2$, $C_{2 \leq 2 \leq 2} = 1$, we get the second relations in \eqref{eqn:L21edges} and \eqref{eqn:L21triangles} from \eqref{eqn:ENRF} and the third and last relations in \eqref{eqn:L21edges} and \eqref{eqn:L21triangles} from \eqref{eqn:ENRS}.
\end{proof}

We write $N_1 = N_1 (\infty)$ and so on for the moments with unlimited radii.
The formulas in Lemma~\ref{lem:moments} simplify since $\gamma(k,\infty) = \Gamma(k) = (k-1)!$.
The Gamma function for half-integers is $\Gamma (k+\sfrac{1}{2}) = \sqrt{\pi} \cdot {(2k)!} / ({4^k k!})$, so $\Gamma (\sfrac{3}{2}) = \frac{1}{2} \sqrt{\pi}$ and $\Gamma (\sfrac{5}{2}) = \frac{3}{4} \sqrt{\pi}$.
We thus get
\begin{align}
  \Expect{N_1} &= 2n ; \hspace{0.4in}
  \Expect{F_1}  = \sqrt{n} ; \hspace{0.4in}
  \Expect{S_1}  = \tfrac{2}{\pi} ; \hspace{0.4in} 
    \label{eqn:L21edgestoo} \\
  \Expect{N_2} &= n ; \hspace{0.47in}
  \Expect{F_2}  = \tfrac{3}{4} \sqrt{n}; \hspace{0.32in}
  \Expect{S_2}  = \tfrac{2}{\pi} .
    \label{eqn:L21trianglestoo}
\end{align}

\subsection{The Poisson Point Process Near the Boundary}
\label{sec:2.2}

We will also need bounds for the critical edges and triangles whose smallest enclosing circles cross the boundary while their centers lie inside $[0,1]^2$.
To this end, we write $N_1^\partial$, $N_2^\partial$ for their numbers, and $F_1^\partial$, $F_2^\partial$ for the sums of their radii, respectively.
\begin{lemma}
  \label{lem:moments_boundary}
  Let $A$ be a stationary Poisson point process with intensity $n > 0$ in $\Rspace^2$.
  Then the critical edges and triangles of the radius function on $\Delaunay{A}$ whose smallest enclosing circles cross the boundary of $[0,1]^2$ while their centers lie in $[0,1]^2$ satisfy
  \begin{align}
    \Expect{N_1^\partial} &\leq 8 \sqrt{n}; \hspace{0.25in}
    \Expect{F_1^\partial}  \leq \tfrac{16}{\pi}; \hspace{0.25in}
    \Expect{N_2^\partial}  \leq 6 \sqrt{n}; \hspace{0.25in}
    \Expect{F_2^\partial}  \leq \tfrac{16}{\pi}.
      \label{eqn:L22}
  \end{align}
\end{lemma}
\begin{proof}
  To prove the first two inequalities in \eqref{eqn:L22}, we consider all critical edges with midpoints in the horizontal strip of unit height, $\Rspace \times [0,1]$.
  For each unit square in this strip, the expected moments of the critical edges with midpoints in the square and radius at most $r_0$ are given in \eqref{eqn:L21edges}.
  For each $t \in \Rspace$, let $n_1 (t)$ be the number of smallest enclosing circles of the considered edges that cross the vertical line passing through the point $(t,0)$, and similarly write $f_1(t)$ for the sum of radii of these edges.
  By symmetry of the Poisson point process, the expected moments are the same for every $t$, so we can compute them by integration.
  An edge contributes along an interval of length twice the radius of the edge.
  Hence,
  \begin{align}
    \Expect{n_1(0)} 
      &= \lim_{x \to \infty} \tfrac{1}{2x} \Expect{ \int_{t=-x}^x n_1 (t) \diff t }
       = \lim_{x \to \infty} \tfrac{1}{2x} \left( 2x \Expect{2 F_1} \right)
       = 2 \cdot \tfrac{2 \sqrt{n}}{\sqrt{\pi}} \cdot \Gamma(\tfrac{3}{2})
       = 2 \sqrt{n} ; 
         \label{eqn:Lemma22N} \\
    \Expect{f_1(0)} 
      &= \lim_{x \to \infty} \tfrac{1}{2x} \Expect{ \int_{t=-x}^x f_1 (t) \diff t }
       = \lim_{x \to \infty} \tfrac{1}{2x} \left( 2x \Expect{2 S_1} \right)
       = 2 \cdot \tfrac{2}{\pi} \cdot \Gamma(2)
       = \tfrac{4}{\pi} ,
         \label{eqn:Lemma22F}
  \end{align}
  where the second equalities in \eqref{eqn:Lemma22N} and \eqref{eqn:Lemma22F} hold because the differences between the integrals and the expectations (which are caused by the finiteness of the strip) divided by $2x$ go to $0$.
  Indeed, both sides accumulate the radii or squared radii of those enclosing circles that fully fit inside $[-x,x] \times \Rspace$, while ignoring the circles that do not intersect the strip.
  The difference is bounded from above by the sum of the radii or squared radii of the remaining enclosing circles, which cross the vertical lines at $-x$ and $x$, i.e., by $4 \Expect{f_1(0)}$ or $4 \Expect{{s_1(0)}}$. 
  The latter are finite numbers, so the limits vanish.

  \smallskip
  Returning to the unit square, this implies that $\Expect{n_1 (0)}$ and $\Expect{f_1 (0)}$ are upper bounds on the expected number and sum of radii of the critical edges whose smallest enclosing circles cross the left side of $[0,1]^2$ and have their centers inside $[0,1]^2$.
  There are four sides of $[0,1]^2$, which implies
  \begin{align}
    \Expect{N_1^\partial} &\leq 4 \Expect{n_1 (0)} = 8 \Expect{F_1} = 8 \sqrt{n}; \hspace{0.3in} \Expect{F_1^\partial} \leq 4 \Expect{f_1 (0)} = 8 \Expect{S_1} = \tfrac{16}{\pi} ;
      \label{eqn:edges_near_boundary} \\
    \Expect{N_2^\partial} &\leq 4 \Expect{n_2 (0)} = 8 \Expect{F_2} = 6 \sqrt{n}; \hspace{0.3in} \Expect{F_2^\partial} \leq 4 \Expect{f_2 (0)} = 8 \Expect{S_2} = \tfrac{16}{\pi} ,
      \label{eqn:triangles_near_boundary}
  \end{align}
  in which we get \eqref{eqn:triangles_near_boundary} by defining $n_2(t)$, $f_2(t)$ analogously to $n_1(t)$, $f_1(t)$ for the critical triangles and applying the same argument to them.
\end{proof}

\subsection{The Uniform Distribution Near the Boundary}
\label{sec:2.3}

We also need a minor generalization of Theorem~1 in \cite{ENR17}, in which we select an edge or triangle if its smallest enclosing circle is empty of points with probability $0 < \eta \leq 1$.
The details of how the proof of Theorem~1 is modified are in \cite[Appendix~A]{DERS25}, and there are simple formulas connecting the new moments to those in \eqref{eqn:ENRN}, \eqref{eqn:ENRF}, and \eqref{eqn:ENRS}:
\begin{align}
  \Expect{N_{\ell \leq k \leq d}^{(\eta)} (r_0)}
    &= \frac{1}{\eta^k} \cdot \frac{\gamma(k,\eta x)}{\gamma(k,x)} \cdot \Expect{N_{\ell \leq k \leq d} (r_0)} ; \\
  \Expect{F_{\ell \leq k \leq d}^{(\eta)} (r_0)}
    &= \frac{1}{\eta^{k+\frac{1}{d}}} \cdot \frac{\gamma(k+\frac{1}{d},\eta x)}{\gamma(k+\frac{1}{d},x)} \cdot \Expect{F_{\ell \leq k \leq d} (r_0)} ; \\
  \Expect{S_{\ell \leq k \leq d}^{(\eta)} (r_0)}
    &= \frac{1}{\eta^{k+\frac{2}{d}}} \cdot \frac{\gamma(k+\frac{2}{d},\eta x)}{\gamma(k+\frac{2}{d},x)} \cdot \Expect{S_{\ell \leq k \leq d} (r_0)} ,
\end{align}
in which $x = \varrho \nu_d r_0^d$, as before.
The incomplete gamma functions cancel each other at $r_0 = \infty$.
Furthermore setting $d=2$, $\varrho = n$, $\eta = \frac{1}{2}$, and $\ell = k = 1$ as well as $\ell = k = 2$, we get
\begin{align}
  \Expect{N_1^{(1/2)}} &= 4n ;
    \hspace{0.4in}
  \Expect{F_1^{(1/2)}}  = \sqrt{8n} ;
    \hspace{0.4in}
  \Expect{S_1^{(1/2)}}  = \tfrac{8}{\pi} ; \\
  \Expect{N_2^{(1/2)}} &= 4n ;
    \hspace{0.4in}
  \Expect{F_2^{(1/2)}}  = \sqrt{18 n};
    \hspace{0.33in}
  \Expect{S_2^{(1/2)}}  = \tfrac{16}{\pi} ;
\end{align}
compare with \eqref{eqn:L21edgestoo} and \eqref{eqn:L21trianglestoo}.
Applying the argument used to prove Lemma~\ref{lem:moments_boundary}, we thus get bounds on the expected $0$-th and $1$-st moments of the critical edges and acute triangles whose smallest enclosing circles are empty with half the probability than before.
This implies bounds for the critical edges and triangles of the Delaunay mosaic of the points inside $[0,1]^2$ near the boundary of the unit square.
Generalizing the earlier notation, we write $N_1^{(1/2) \partial}$ for the number of critical edges whose smallest enclosing circles are empty with probability $\frac{1}{2}$ and cross the boundary of $[0,1]^2$, while its midpoint lies in $[0,1]^2$, etc.
\begin{lemma}
  \label{lem:moments_boundary_too}
  Let $A$ be a stationary Poisson point process with intensity $n > 0$ in $\Rspace^2$.
  Then the critical edges and triangles of the radius function on $\Delaunay{A \cap [0,1]^2}$ whose smallest enclosing circles cross the boundary of $[0,1]^2$ satisfy
  \begin{align}
    \Expect{N_1^{(1/2) \partial}} &\leq \sqrt{512n};
      \hspace{0.18in}
    \Expect{F_1^{(1/2) \partial}} \leq \tfrac{64}{\pi};
      \hspace{0.18in}
    \Expect{N_2^{(1/2) \partial}} \leq \sqrt{1152 n};
      \hspace{0.18in}
    \Expect{F_2^{(1/2) \partial}} \leq \tfrac{128}{\pi}.
      \label{eqn:L23}
  \end{align}
\end{lemma}
\begin{proof}
  A critical edge in the Delaunay mosaic of $A \cap [0,1]^2$ has two antipodal points on the smallest enclosing circle, which implies that at least half the area of the disk bounded by this circle belongs to $[0,1]^2$.
  By definition, this disk does not contain any points of $A \cap [0,1]^2$.
  It may however contain points of $A$, namely if they lie outside $[0,1]^2$.
  It follows that the probability of this circle being empty of points in $A \cap [0,1]^2$ is at least twice the probability that it is empty of points in $A$.
  In analogy to the proof of Lemma~\ref{lem:moments_boundary}, we therefore get
  \begin{align}
    \Expect{N_1^{(1/2) \partial}} &\leq 8 \Expect{F_1^{(1/2)}} = \sqrt{512 n} ; \hspace{0.37in}
    \Expect{F_1^{(1/2) \partial}} \leq 8 \Expect{S_1^{(1/2)}} = \tfrac{64}{\pi} ;
    \\
    \Expect{N_2^{(1/2) \partial}} &\leq 8 \Expect{F_2^{(1/2)}} = \sqrt{1152 n} ; \hspace{0.30in}
    \Expect{F_2^{(1/2) \partial}} \leq 8 \Expect{S_2^{(1/2)}} = \tfrac{128}{\pi} ,
  \end{align}
  where we get the second two relations by applying the same argument to the critical triangles of the radius function on $\Delaunay{A \cap [0,1]^2}$.
\end{proof}

\section{Euclidean Minimum Spanning Trees}
\label{sec:3}

This section presents the main result of this paper, which is a new lower bound for the asymptotic constant of the expected length of the Euclidean minimum spanning tree in the plane.
We begin with a brief introduction of the relevant concepts and a review of the prior work on this constant.

\subsection{Expected Length}
\label{sec:3.1}

For a finite set, $A \subseteq \Rspace^2$, the \emph{Euclidean minimum spanning tree}, or \emph{EMST} for short, is the tree that connects all points and minimizes the sum of edge lengths.
It is not necessarily unique, but its length is.
Kruskal's algorithm \cite{Kru56} constructs this tree by considering the edges with endpoints in $A$ in the order of their lengths, and adds an edge to the tree if it connects two yet different connected components.
Letting $n$ be the number of points in $A$, the EMST can be constructed in $O(n \log n)$ time by first computing the Delaunay mosaic of $A$ and second running Kruskal's algorithm on the edges of this mosaic.
It is indeed not difficult to see that the edges in the tree all belong to this mosaic and, assuming the generic case in which all edges have distinct lengths, they are in fact the death-giving critical edges of its radius function.

\smallskip
In this paper, we focus on the length of the EMST.
More than half a century ago Beardwood, Hammersley, and Halton~\cite{BHH59} proved that there exists a constant, $\Cmst$, such that the expected length of the EMST of $n$ random points in $[0,1]^2$ is $\Cmst \sqrt{n}$ in the limit, when $n$ goes to infinity.
We refer to $\Cmst$ as the \emph{asymptotic constant} of the expected length of the EMST in $\Rspace^2$.
The notation for this constant is $\Cmst = c(1,2) = \beta (2)$ in \cite{Ste88} and \cite{BHH59}, respectively, with the parameters being the power of the Euclidean edge length that is minimized and the dimension of the ambient space.
We focus on the length (with power $1$) in the Euclidean plane and therefore drop these parameters from the notation.

\smallskip
Although we know that the asymptotic constant exists, its exact value is still unknown.
As mentioned in the introduction, the best upper and lower bounds prior to this paper are $0.6008 \leq \Cmst \leq 0.7072$.
More precisely, the upper bound by Gilbert~\cite{Gil65} is $\sfrac{\sqrt{2}}{2}$ and the lower bound by Avram and Bertsimas~\cite{AvBe92} is extracted from a formula to calculate $\Cmst$, which however is a series expansion whose terms require difficult numerical integrations.
We will improve the lower bound with a relatively elementary argument shortly.
There are three sources in the literature that estimate the asymptotic constant with computational experiments, suggesting $\Cmst$ is approximately $0.6331$, $0.656$, and $0.68$ by Cortina-Borja and Robinson \cite{CBRo00}, Roberts \cite{Rob68}, and Gilbert \cite{Gil65}, respectively.
All these estimates are larger than our lower bound; see Theorem~\ref{thm:improvement}.
The results of our computational experiments described in Section~\ref{sec:6} suggest that $\Cmst$ is between $0.646$ and $0.648$ and thus lies between the two smaller of the three estimates in the literature.

\subsection{Improved Lower Bound for Asymptotic Constant}
\label{sec:3.2}

Before stating and proving the main result of this section, we recall a connection between a stationary Poisson point process with intensity $n$ in $[0,1]^2$ and a set of $n$ points chosen uniformly at random in $[0,1]^2$.
In contrast to the uniform distribution model, the Poisson point process does not necessarily have $n$ points but rather $k \geq 0$ points with probability $\Prob{k} = n^k e^{-n} / k!$, and therefore $\sum_{k \geq 0} k \cdot \Prob{k} = n$ points in expectation.
Suppose $f(n)$ is the expectation of a random variable---e.g., the length of the EMST---of $n$ points chosen uniformly at random in $[0,1]^2$.
Is this also the expectation of that random variable for the stationary Poisson point process with intensity $n$ in $[0,1]^2$, expect for a lower order term?
In other words, we ask whether
\begin{align}
  \sum\nolimits_{k \geq 0} f(k) \cdot \Prob{k} &= f(n) + o(f(n)).
    \label{eq:poisson-vs-uniform}
\end{align}
The general results in \cite[Theorem~2.1 and Section~4]{CGKK13} imply that \eqref{eq:poisson-vs-uniform} holds for $f(n) = \Cmst \sqrt{n}$.
In the case of the EMST length, we know that there exists a constant, $\Cmst$, such that the expectation is $f(n) = c \sqrt{n}$ in the limit. 
Then \eqref{eq:poisson-vs-uniform} also holds in the limit: if $n$ is large, we can approximate $f(k)$ by $c\sqrt{k}$ for all $k>n/2$, and disregard the entries for smaller $k$ as their sum is overall of order less than $\sqrt{n}$.
Therefore, when we prove in Theorem~\ref{thm:improvement} that the expectation of the left-hand side of \eqref{eq:poisson-vs-uniform} is bounded from below by $c_0 \sqrt{n}$ in the limit, then we have $c_0 \sqrt{n} \leq c \sqrt{n} + o(\sqrt{n})$ for large $n$, and $c_0 \leq \Cmst$ follows.

\smallskip
Let us focus on the random variable that maps a finite set of points to the length of the Euclidean minimum spanning tree.
We know that the EMST of $n+1$ points consists of $n$ edges, they all belong to the Delaunay mosaic and, more specifically, they are critical edges of the radius function on the Delaunay mosaic of the points.
Furthermore, Lemma~\ref{lem:moments} spells out how long the critical edges are in expectation.
Ignoring for now that this lemma talks about a Poisson point process rather than points chosen uniformly at random in $[0,1]^2$, we can just take the $n$ shortest critical edges and use their expected length to glean a lower bound for the expected length of the EMST.
Slightly less directly, we find the threshold, $r_0$, for which we get an expected $n$ critical edges with length at most $2 r_0$, and compute $F_1(r_0)$.
To see that $2 \Expect{F_1(r_0)}$ is a lower bound for the expected length of the $n$ shortest edges, we note that we over-count as many edges as we under-count, in expectation.
By construction, the over-counted edges have length at most $2 r_0$, while the under-counted edges are longer than $2 r_0$, so the result of the computation cannot be more than the expected length of the $n$ shortest edges.
Following this strategy, we set $r_0 = \sqrt{\ln 2 /(n \pi)}$ and $x = n \pi r_0^2 = \ln 2$, so that the expected number of critical edges that have length at most $2 r_0$, and the expected total length of these edges are
\begin{align}
  \Expect{N_1(r_0)}
      &= 2n \cdot \gamma (1, x) 
       = 2n \int\nolimits_{t=0}^x e^{-t} \diff t
       = -2n \left[ e^{-t} \right]_0^x
       = 2n - 2n e^{-x} 
       = n ; 
    \label{eqn:envelope1} \\
  \Expect{F_1 (r_0)}
      &= \tfrac{2 \sqrt{n}}{\sqrt{\pi}} \cdot \gamma(\tfrac{3}{2}, x)
       = \tfrac{2 \sqrt{n}}{\sqrt{\pi}} \int\nolimits_{t=0}^x t^{1/2} e^{-t} \diff t
       = 0.2912\ldots \cdot \sqrt{n} ,
    \label{eqn:envelope2}
\end{align}
in which we resort to numerical integration to the get the right-hand side of \eqref{eqn:envelope2}.
In words, for the threshold that gives $n$ edges in expectation, we get an expected length equal to twice $0.2912\ldots \cdot \sqrt{n}$, but twice $0.2912$ is less than $0.6008$, so making this argument rigorous is futile.
However, among the shortest edges, there are many that form cycles and can therefore not serve as edges of the EMST.
These edges need to be replaced by further edges in the pool of candidates to be considered for the EMST.
Hence, this tree must be longer than what we just estimated.
The main new idea is to use the critical triangles of the radius function on the Delaunay mosaic to get a handle on the edges that form cycles.
This is possible because every critical triangle implies the existence of a corresponding critical edge with smaller radius that gives birth to the loop the triangle fills.
We now make these ideas concrete.
\begin{theorem}
  \label{thm:improvement}
  Let $\Cmst$ be the asymptotic constant of the planar Euclidean minimum spanning tree; that is: the constant such that $\Cmst \sqrt{n}$ is the expected length of the EMST of $n$ points sampled uniformly at random from $[0,1]^2$, in the limit, when $n$ goes to infinity.
  Then $0.6289\ldots \leq \Cmst$.
\end{theorem}
\begin{proof}
  Rather than $n$ points chosen uniformly at random in $[0,1]^2$, we consider a stationary Poisson point process with intensity $n$ in $\Rspace^2$, observe that the subset of points in $[0,1]^2$ is a stationary Poisson point process with the same intensity in $[0,1]^2$, and refer to \eqref{eq:poisson-vs-uniform} to justify this choice.
  The main new insight differentiating this proof from the earlier back-of-the-envelope argument is that a birth-giving critical edge is identified as such by the paired critical triangle.
  If these edges are short, then their exclusion from the pool of possible edges forces the use of longer edges in the minimum spanning tree.
  To implement this idea, let $r_0 = {1}/{\sqrt{n \pi}}$ and $x = n \pi r_0^2 = 1$.
  According to Lemma~\ref{lem:moments}, the critical edges and triangles of radius at most $r_0$ satisfy
  \begin{align}
    \begin{split}
    \Expect{N_1 (r_0)} - \Expect{N_2 (r_0)}
      &= 2n \int\nolimits_{t=0}^x e^{-t} \diff t - n \int\nolimits_{t=0}^x t e^{-t} \diff t \\
      &= 2n (1-e^{-x})  - n (1- (x+1) e^{-x})
       = n ; 
    \end{split}
    \label{eqn:expN} \\
    \Expect{F_1 (r_0)} - \Expect{F_2 (r_0)}
      &= \tfrac{2 \sqrt{n}}{\sqrt{\pi}} \int\nolimits_{t=0}^x t^{1/2} e^{-t} \diff t - \tfrac{\sqrt{n}}{\sqrt{\pi}} \int\nolimits_{t=0}^x t^{3/2} e^{-t} \diff t 
       = 0.31445\ldots \cdot \sqrt{n} ,
    \label{eqn:expF}
  \end{align}
  in which we get the right-hand side in \eqref{eqn:expN} by noticing $2e^{-x} = (x+1) e^{-x}$, and the right-hand side in \eqref{eqn:expF} by numerical integration.
  In words, for this threshold we expect a surplus of $n$ critical edges beyond the expected number of critical triangles, and a surplus of $0.6289\ldots \cdot \sqrt{n}$ in length.
  By selection, the radii of these critical triangles are at most $r_0$, so the corresponding birth-giving edges are indeed short.

  \smallskip
  This is promising but we need to cope with the difference between the Poisson point process, $A \subseteq \Rspace^2$, and its subset, $U = A \cap [0,1]^2$.
  In particular, \eqref{eqn:expN} and \eqref{eqn:expF} apply to the critical edges and critical triangles of the radius function on $\Delaunay{A}$ whose centers lie inside $[0,1]^2$, which are not necessarily the critical edges and critical triangles of the radius function on $\Delaunay{U}$.
  Let $E_A, T_A$ and $E_U, T_U$ denote the sets of such critical edges and triangles in $\Delaunay{A}$ and $\Delaunay{U}$, respectively, and note that the differences between these sets are concentrated near the boundary of $[0,1]^2$.
  The main argument will be about the critical edges and triangles shared by the two mosaics.

  \smallskip
  To make this concrete, we use the classification of $E_U$ into birth-giving and death-giving edges, and the pairing of the former with the triangles in $T_U$ provided by the persistent homology of the radius function on $\Delaunay{U}$; see e.g.\ \cite[Chapter~VII]{EdHa10}.
  We write $E \subseteq E_U$ and $T \subseteq T_U$ for the subsets obtained by excluding edges and triangles that are not also critical in $\Delaunay{A}$ or are paired with triangles and edges not in $\Delaunay{A}$.
  In other words, $E = (E_U \cap E_A) \setminus E''$ and $T = (T_U \cap T_A) \setminus T''$, in which the edges in $E'' \subseteq E_U$ are paired with triangles in $T'= T_U \setminus T_A$ and the triangles in $T'' \subseteq T_U$ are paired with edges in $E' = E_U \setminus E_A$.
  Finally, set $E''' = E_A \setminus E_U$ and $T''' = T_A \setminus T_U$.
  Note that edges and triangles in $E$ and $T$ have the same radii in $\Delaunay{U}$ and $\Delaunay{A}$.
  To bound the numbers and sums of radii, we observe that all edges and triangles in $E_U$ and $T_U$ whose smallest enclosing circles are contained in $[0,1]^2$ also belong to $E_A$ and $T_A$, respectively.
  Hence,
  \begin{align}
    \Expect{N_{E'}} &\leq \sqrt{512 n};
      \hspace{0.3in}
    \Expect{F_{E'}}  \leq \tfrac{64}{\pi};
      \hspace{0.3in}
    \Expect{N_{T'}}  \leq \sqrt{1152 n};
      \hspace{0.27in}
    \Expect{F_{T'}}  \leq \tfrac{128}{\pi};
      \label{eqn:ETprime} \\
    \Expect{N_{E''}} &\leq \sqrt{1152 n};
      \hspace{0.21in}
    \Expect{F_{E''}}  \leq \tfrac{128}{\pi};
      \hspace{0.21in}
    \Expect{N_{T''}}  \leq \sqrt{512 n};
      \label{eqn:ETprimeprime} \\
    \Expect{N_{E'''}} &\leq \sqrt{64 n};
      \hspace{0.32in}
    \Expect{F_{E'''}}  \leq \tfrac{16}{\pi};
      \hspace{0.24in}
    \Expect{N_{T'''}}  \leq \sqrt{36 n};
      \hspace{0.33in}
    \Expect{F_{T'''}}  \leq \tfrac{16}{\pi},
      \label{eqn:ETprimeprimeprime}
  \end{align}
  in which \eqref{eqn:ETprime} is implied by Lemma~\ref{lem:moments_boundary_too} and \eqref{eqn:ETprimeprimeprime} is implied by Lemma~\ref{lem:moments_boundary}.
  The relations in \eqref{eqn:ETprimeprime} are implied by those in \eqref{eqn:ETprime}.
  In particular, the second inequality in \eqref{eqn:ETprimeprime} is implied by the last inequality in \eqref{eqn:ETprime} because the radius of a birth-giving edge is at most the radius of the death-giving triangle the edge is paired with.
  Since all cardinalities are at most $O (\sqrt{n})$ and all sums of radii are at most $O (1)$, in expectation, discarding these edges and triangles barely changes anything, as we will see shortly.
  Next, we partition $E$ into five sets:
  \begin{itemize}
    \item \emph{short}: $\Dshort$ are the death-giving edges whose radii are at most $r_0$;
    \item \emph{long}: $\Dlong$ are the death-giving edges whose radii exceed $r_0$;
    \item \emph{certified short}: $\Bcertified$ are the birth-giving edges paired with triangles whose radii are at most $r_0$;
    \item \emph{uncertified short}: $\Buncertified$ are the birth-giving edges whose radii at most $r_0$ and are paired with triangles whose radii exceed $r_0$.
  \end{itemize}
  There are still the \emph{long} birth-giving edges in $E$---whose radii exceed $r_0$---but we will not need them.
  Writing $N_{\Dshort}$ for the number of edges in $\Dshort$ and $F_{\Dshort}$ for the sum of their radii, etc., we get
  \begin{align}
    N_1 (r_0) - N_{E''} - N_{E'''} &\leq N_{\Dshort} + N_{\Bcertified} + N_{\Buncertified} 
        \leq N_1 (r_0) ;
      \label{eqn:NDBBprime} \\
    N_2 (r_0) - N_{T''} - N_{T'''} &\leq \hspace{0.48in} N_{\Bcertified} \hspace{0.48in} \leq N_2 (r_0); 
      \label{eqn:NB} \\
    F_1 (r_0) - F_{E''} - F_{E'''} &\leq F_{\Dshort} + \hspace{0.035in} F_{\Bcertified} + \hspace{0.03in} F_{\Buncertified} \hspace{0.01in} \leq F_1 (r_0); 
      \label{eqn:FDBBprime} \\
    &{} \hspace{0.68in}  F_{\Bcertified} \hspace{0.50in} \leq F_2 (r_0) .
      \label{eqn:FB}
  \end{align}
  Here \eqref{eqn:NDBBprime} says that of the critical edges in $\Delaunay{A}$ with radius at most $r_0$ whose centers lie in $[0,1]^2$, the sets $\Dshort$, $\Bcertified$, $\Buncertified$ miss only the edges in $E''$ and $E'''$, of which there are at most $O(\sqrt{n})$ in expectation by \eqref{eqn:ETprimeprime} and \eqref{eqn:ETprimeprimeprime}.
  Similarly, \eqref{eqn:NB} says that of all critical edges paired with triangles in $\Delaunay{A}$ with radius at most $r_0$, the set $\Bcertified$ is paired with all triangles in $T_A$ except the ones in $T''$ and $T'''$, of which there are at most $O(\sqrt{n})$ in expectation by \eqref{eqn:ETprimeprime} and \eqref{eqn:ETprimeprimeprime}.
  Relation \eqref{eqn:FDBBprime} follows because \eqref{eqn:NDBBprime} can also be written in terms of sets and inclusions of sets, and relation \eqref{eqn:FB} follows from \eqref{eqn:NB} and the fact that the radius of an edge is at most the radius of the paired triangle.
  Subtracting \eqref{eqn:NB} from \eqref{eqn:NDBBprime} and \eqref{eqn:FB} from \eqref{eqn:FDBBprime}, we get
  \begin{align}
    N_1 (r_0) - N_2 (r_0) - N_{E''} - N_{E'''} &\leq N_{\Dshort} + N_{\Buncertified} \leq N_1 (r_0) - N_2 (r_0) + N_{T''} + N_{T'''} ; 
      \label{eqn:NDBprime} \\
    F_1 (r_0) - F_2 (r_0) - F_{E''} - F_{E'''} &\leq F_{\Dshort} + F_{\Buncertified} ,
      \label{eqn:FDBprime}
  \end{align}
  in which all primed terms are small in expectation, so that their presence in the above inequalities makes little difference.
  Therefore after taking expectations, $\Dshort$ and $\Buncertified$ have about the right number of edges, namely $n \pm O (\sqrt{n})$, as computed in \eqref{eqn:expN}.
  Similarly, $\Dshort$ and $\Dlong$ have $n \pm O (\sqrt{n})$ edges in expectation, which implies that $\Buncertified$ and $\Dlong$ have about the same number of edges in expectation.
  By definition, the edges in $\Dlong$ are longer than the ones in $\Buncertified$, but $\Buncertified$ may have more edges than $\Dlong$, so an inequality for the sum of radii is not immediate.
  However, the difference in numbers is at most $O (\sqrt{n})$, so we get
  \begin{align}
    \Expect{F_{\Buncertified}} &\leq \Expect{F_{\Dlong}} \pm r_0 O (\sqrt{n})
      = \Expect{F_{\Dlong}} \pm O (1)
  \end{align}
  because $r_0 \sqrt{n} = 1 / \sqrt{\pi}$ is only a constant.
  Rewriting \eqref{eqn:FDBprime} with expectations thus gives
  \begin{align}
    \Expect{F_1 (r_0)} - \Expect{F_2 (r_0)} \pm O (1)
      &\leq \Expect{F_{\Dshort}} + \Expect{F_{\Buncertified}}
       \leq \Expect{F_{\Dshort}} + \Expect{F_{\Dlong}} \pm O (1) ,
  \end{align}
  in which we use \eqref{eqn:ETprimeprime} and \eqref{eqn:ETprimeprimeprime} to bound the expectations of $F_{E''}$ and $F_{E'''}$.
  By \eqref{eqn:expF}, the left-hand side is $0.31445\ldots \cdot \sqrt{n} \pm O (1)$.
  The claimed inequality for the constant $\Cmst$ follows because $2F_{\Dshort} + 2F_{\Dlong}$ is a lower bound for the length of the EMST of $U$, and the asymptotic constant is the same no matter whether the expectation is over $n$ points sampled uniformly in $[0,1]^2$ or over a stationary Poisson point process with intensity $n$ in $[0,1]^2$, as explained at the beginning of this subsection.
\end{proof}

\section{Points with Two Colors}
\label{sec:4}

This section presents the needed background in persistent homology and its chromatic variant, and we refer to \cite{EdHa10} and \cite{CDES23} for more comprehensive treatments of these topics.
In addition, we introduce a novel generalization of the EMST.

\subsection{Persistent Homology}
\label{sec:4.1}

A \emph{filtration} is a linear sequence of topological spaces, in which each space is a subspace of the next.
To construct an example, we may start with a complex, $K_0$, for $i > 0$ get $K_i$ from $K_{i-1}$ by adding a cell whose proper faces are already in $K_{i-1}$, and end with $K_n$.
Applying the $p$-th homology functor for field coefficients in $\kkk$, we get a sequence of vector spaces connected from left to right by linear maps, $\Hgroup{p} (K_0) \to \Hgroup{p} (K_1) \to \ldots \to \Hgroup{p} (K_n)$, which is sometimes referred to as a \emph{persistence module}.
Importantly, it decomposes into \emph{interval modules} of the form $\ldots \to 0 \to \kkk \to \ldots \to \kkk \to 0 \to \ldots$, in which the maps connecting the copies of $\kkk$ are identities and all others are zero maps.
Letting $i$ and $j-1$ be the positions of the first and last copy of $\kkk$, the interval module represents a degree-$p$ homology class that is \emph{born} at $K_i$ and \emph{dies entering} $K_j$.
We say the cell in $K_i \setminus K_{i-1}$ \emph{gives birth} and the cell in $K_j \setminus K_{j-1}$ \emph{gives death} to this class.

\smallskip
In this paper, the $K_i$ are primarily sublevel sets of the radius function on the Delaunay mosaic of $A \subseteq \Rspace^2$, denoted $f \colon \Delaunay{A} \to \Rspace$.
Let $r_0 < r_1 < \ldots$ be the values of $f$ and recall that $K_i = \Radiusf^{-1} [0,r_i]$ is the alpha complex of $A$ and radius $r_i$.
Assume for convenience that each $K_i$ differs from $K_{i-1}$ by the simplices in a single interval of $f$.
If this interval contains two or more and therefore non-critical simplices, then $K_i$ and $K_{i-1}$ have the same homotopy type and thus isomorphic homology groups in all degrees.
On the other hand, if this interval contains a single and therefore critical simplex, then there is exactly one degree for which the homology groups of $K_i$ and $K_{i-1}$ are not isomorphic:
letting $p$ be the dimension of the critical simplex, then 
either $\rank{\Hgroup{p-1}(K_i)} = \rank{\Hgroup{p-1}(K_{i-1})} - 1$ (the simplex gives death to a $(p-1)$-cycle), or $\rank{\Hgroup{p}(K_i)} = \rank{\Hgroup{p}(K_{i-1})} + 1$ (the simplex gives birth to a $p$-cycle).
For example, a critical edge gives death to a component iff it belongs to the EMST of the points, and it gives birth to a loop iff it does not belong to this tree.
Being critical or non-critical is a local property, while giving birth or death is a global property.
The methods from probability theory are well suited to study this local property---see for example Lemma~\ref{lem:moments}, which gives the expected densities for a Poisson point process---while the global property is probabilistically elusive and only incomplete information is available---see the experimental findings documented in \cite{ENOS20}.
Indeed the authors of this paper believe that the difficulty in determining the constant $\Cmst$ is tied up with the difference between birth-giving and death-giving critical edges of the radius function on the Delaunay mosaic.

\subsection{Chromatic Persistent Homology}
\label{sec:4.2}

A recent development in the theory of persistent homology uses colors to distinguish different types of points.
We limit ourselves to the case of two colors, $0$ and $1$, and write $\chi \colon A \to \{0,1\}$ for the $2$-coloring.
A useful property is the following: if $A \subseteq \Rspace^2$ is a stationary Poisson point process with intensity $n$, and the points are randomly colored $0$ or $1$, with probability $\tfrac{1}{2}$ each, then $A_0 = \chi^{-1}(0)$ and $A_1 = \chi^{-1}(1)$ are stationary Poisson point processes with intensity $\tfrac{n}{2}$ each; see Kingman~\cite{Kin93} for background on Poisson point processes.

\smallskip
To study the interaction between the colors, write $A(r)$ for the points at distance at most $r$ from at least one point in $A$.
Equivalently, $A(r)$ is the union of closed disks of radius $r$ centered at the points in $A$.
We extend this notation to the mono-chromatic sets, $A_0$ and $A_1$, and observe that $A_0(r) \subseteq A(r)$ and also $A_1(r) \subseteq A(r)$.
Taking the disjoint union of the two mono-chromatic sets, we get two parallel filtrations with vertical inclusion maps connecting the sets for the same radius.
After applying the homology functor for field coefficients, we get a commutative ladder of vector spaces for each degree:
\begin{align}
  \begin{array}{ccccccccccc}
    \ldots & \to & \Hgroup{p}(A(r)) & \to & \Hgroup{p}(A(s)) & \to & \ldots & \to & \Hgroup{p}(A(t)) & \to & \ldots \\
    & & \uparrow & & \uparrow & & & & \uparrow & & \\
    \ldots & \to & \!\!\!\Hgroup{p}(A_0(r) \sqcup A_1(r))\!\!\! & \to & \!\!\!\Hgroup{p}(A_0(s) \sqcup A_1(s))\!\!\! & \to & \ldots & \to & \!\!\!\Hgroup{p}(A_0(t) \sqcup A_1(t))\!\!\! & \to & \ldots
  \end{array}
\end{align}
Letting $\kappa_r \colon \Hgroup{p}(A_0(r) \sqcup A_1(r)) \to \Hgroup{p}(A(r))$ be the linear map induced by the inclusions, we get six (1-parameter) persistence modules, one each for the domains, codomains, pairs, images, cokernels, and kernels.
For each integer $p$, we write $\Domain{p}{f}$, $\Codomain{p}{f}$, $\Relative{p}{f}$, $\Image{p}{f}$, $\Cokernel{p}{f}$, and $\Kernel{p}{f}$ for the corresponding \emph{degree-$p$ chromatic persistence diagrams}, and $\norm{\Domain{p}{f}}_1$ for the $1$-norm of the degree-$p$ domain diagram, etc.
The short exact sequence that connects the homology of the kernel, domain, image implies a linear relation for the $1$-norms of the corresponding persistence diagrams, and similarly for the other groups:
\begin{align}
  \norm{\Kernel{p}{f}}_1 + \norm{\Image{p}{f}}_1
    &= \norm{\Domain{p}{f}}_1 ;
  \label{eqn:short1} \\
  \norm{\Image{p}{f}}_1 + \norm{\Cokernel{p}{f}}_1
    &= \norm{\Codomain{p}{f}}_1 ;
  \label{eqn:short2} \\
  \norm{\Cokernel{p}{f}}_1 + \norm{\Kernel{p-1}{f}}_1
    &= \norm{\Relative{p}{f}}_1,
  \label{eqn:short3}
\end{align}
There are three relations for each degree, $p$, but note the shift to degree $p-1$ for the kernel in \eqref{eqn:short3}.
This shift is a consequence of the long exact sequence of a pair, which lists the homology groups in a spiral rather than a cyclic sequence; see \cite{CDES23} for further details.

\subsection{The Lunar EMST}
\label{sec:4.3}

Of the six chromatic persistence diagrams, the one capturing the persistent relative homology of the pair $(A(r), A_0(r) \sqcup A_1(r))$ is particularly interesting for the analysis of the $1$-norms. 
We reformulate the degree-$1$ diagram and its $1$-norm in more elementary terms.
In our setting, the sets $A_0$ and $A_1$ are required to be disjoint, but this condition is not needed and could be dropped.

\smallskip
Letting $a \in A_0$ and $b \in A_1$, we call the points whose distances from $a$ and $b$ are both at most $r$ the \emph{lune} of radius $r$.
It is the intersection of the two closed disks of radius $r$ centered at the two points, which is non-empty iff $r$ is at least half the distance between $a$ and $b$.
Given $A_0, A_1 \subseteq \Rspace^2$ and $r \geq 0$, we take the union of such lunes over all pairs $(a,b) \in A_0 \times A_1$.
Increasing $r$ from $0$ to $\infty$, we get a filtration and therefore a persistence diagram.
Just like the EMST captures the sequence of merge events between the components of a union of growing disks, we introduce the \emph{lunar EMST} to capture the merge events between the components of the union of growing lunes.
An important difference to the disks is that the lunes start out empty and become non-empty only later during the process.
When a lune wakes up, it may or may not start a new component in the union of lunes.
To quantify the process, we assign twice the radius at which a component is born (its first lune wakes up) to the corresponding vertex, and twice the radius at which two components merge to the corresponding edge.
The factor $2$ is introduced to represent the distance rather than half of it, which is the radius.
The \emph{cost} of the lunar EMST is then the sum of edge costs minus the sum of vertex costs plus the minimum vertex cost, and because we use the factor $2$, this cost is the length of the standard EMST if $A_0 = A_1$
  
\smallskip
We refer to \cite{DERS25} for the details of this construction, and for the proof that the cost of the lunar EMST is indeed twice the $1$-norm of the degree-$1$ chromatic persistence diagram of the pair $(A(r), A_0(r) \sqcup A_1(r))$.
Importantly, the same paper also shows that the expected cost of the lunar EMST of $n$ points chosen uniformly at random in $[0,1]^2$ and randomly $2$-colored converges to a constant times $\sqrt{n}$.
We call this the \emph{asymptotic constant} of the expected cost, denoted $\CLmst$.
Similar to the case of the standard EMST, $\CLmst$ is known to exist but its value is not known.

\section{Expected 1-norms}
\label{sec:5}

The focus of this paper is the expected $1$-norms of the chromatic persistence diagrams.
For the setting at hand, they depend on the expected length of the EMST and the expected cost of the lunar EMST.

\subsection{The Asymptotic Constants}
\label{sec:5.1}

As introduced in Section~\ref{sec:4.2}, we have six chromatic persistence diagrams for the map from the disjoint union of the two mono-chromatic sets to the bi-chromatic set.
We further distinguish between different degrees, which results in eleven diagrams with possibly non-zero $1$-norms.
We will see shortly that the expected $1$-norms are again of the form a constant times $\sqrt{n}$, in the limit, in which six of these constants depend solely on $\Cmst$.
The remaining five expected $1$-norms depend on $\CLmst$, or on $\Cmst$ and $\CLmst$; see Table~\ref{tab:6-pack}.
{\renewcommand{\arraystretch}{1.3}
\begin{table}[htb]
  \centering  \footnotesize
  \begin{tabular}{l||l|l|l}
    & \multicolumn{1}{c|}{$\Kernel{}{f}$} & \multicolumn{1}{c|}{$\Relative{}{f}$} & \multicolumn{1}{c}{$\Cokernel{}{f}$} \\ \hline \hline
    $p=0$ & $\Ckernel{0} = \tfrac{1}{2} (\sqrt{2}-1) \Cmst$ & $\Crelative{0} = 0$ & $\Ccokernel{0} = 0$ \\
    $p=1$ & $\Ckernel{1} = \tfrac{1}{2} \CLmst -\tfrac{1}{4} (\sqrt{2}-1)$ & $\Crelative{1} = \tfrac{1}{2} \CLmst$ & $\Ccokernel{1} = \tfrac{1}{2} \CLmst - \tfrac{1}{2} (\sqrt{2}-1) \Cmst$ \\
    $p=2$ & $\Ckernel{2} = 0$ & $\Crelative{2} = \tfrac{1}{2} \CLmst - \tfrac{1}{4} (\sqrt{2}-1)$ & $\Ccokernel{2} = 0$ \\ \hline \hline
    & \multicolumn{1}{c|}{$\Domain{}{f}$} & \multicolumn{1}{c|}{$\Image{}{f}$} & \multicolumn{1}{c}{$\Codomain{}{f}$} \\ \hline \hline
    $p=0$ & $\Cdomain{0} = \tfrac{1}{2} \sqrt{2} \Cmst$ & $\Cimage{0} = \tfrac{1}{2} \Cmst$ & $\Ccodomain{0} = \tfrac{1}{2} \Cmst$ \\
    $p=1$ & $\Cdomain{1} = \sqrt{2} (\tfrac{1}{2} \Cmst - \tfrac{1}{4})$ & $\Cimage{1} = \tfrac{1}{2} \sqrt{2} \Cmst - \tfrac{1}{4} - \tfrac{1}{2} \CLmst$ & $\Ccodomain{1} = \tfrac{1}{2} \Cmst - \tfrac{1}{4}$ \\
    $p=2$ & $\Cdomain{2} = 0$ & $\Cimage{2} = 0$ & $\Ccodomain{2} = 0$
  \end{tabular} \vspace{0.0in}
  \caption{\footnotesize The constants for the expected $1$-norms of the chromatic persistence diagrams of randomly $2$-colored points in $[0,1]^2$, in which the two colors are equally likely.}
  \label{tab:6-pack}
\end{table}}
\begin{theorem}
  \label{thm:constants}
  Let $A$ be $n$ points sampled uniformly at random from $[0,1]^2$, $\chi \colon A \to \{0,1\}$ a random $2$-coloring, and $\kappa_r$ the map on homology induced by the inclusions of $A_0(r)$ and $A_1(r)$ into $A(r)$.
  Then for $0 \leq p \leq 2$, the expected $1$-norms of the degree-$p$ chromatic persistence diagrams satisfy
  $$ 
    \begin{array}{ccc}
      \norm{\Kernel{p}{f}}_1 = \Ckernel{p} \sqrt{n},
        & \norm{\Relative{p}{f}}_1 = \Crelative{p} \sqrt{n},
        & \norm{\Cokernel{p}{f}}_1 = \Ccokernel{p} \sqrt{n}, \\
      \norm{\Domain{p}{f}}_1 = \Cdomain{p} \sqrt{n},
        & \norm{\Image{p}{f}}_1 = \Cimage{p} \sqrt{n},
        & \norm{\Codomain{p}{f}}_1 = \Ccodomain{p} \sqrt{n}
    \end{array}
  $$
  in the limit when $n$ goes to infinity, with the constants as shown in Table~\ref{tab:6-pack}.
\end{theorem}
\begin{proof}
  We prove the eleven relations in three steps: first for all degree-$0$ diagrams, second for the degree-$1$ domain and codomain diagrams, and third for the remaining diagrams.
  These relations imply that the $1$-norms are indeed of the form a constant times $\sqrt{n}$.
  Recall that $A_0 (r)$ and $A_1(r)$ are the two mono-chromatic unions of disks, and $A(r) = A_0(r)\cup A_1(r)$ is the bi-chromatic union of disks.
  Recall also that the alpha complexes for this radius inside $\Delaunay{A_0}$, $\Delaunay{A_1}$, and $\Delaunay{A}$ have the same homotopy type as the corresponding unions of disks.

  \smallskip \noindent
  \textsc{Degree-$0$ diagrams.}
  Since $A = A_0 \cup A_1$, the cokernel and relative diagrams in degree $0$ are necessarily empty, so their $1$-norms vanish.
  The other four degree-$0$ diagrams are related to the EMST of the points, as we now describe, but see \cite[Appendix~A]{CDES24} for more detail.
  Kruskal's algorithm constructs the EMST of a finite set of points by incrementally adding the shortest edge that connects two yet separate connected components.
  In the language of persistent homology of the growing space, $A(r)$, half the length of this edge is the radius at which one of its components dies.
  Hence, the sum of radii of the death-giving edges is half the length of the EMST.
  Since all components are born at radius $r = 0$, this implies that the $1$-norm of the degree-$0$ persistence diagram is half the length of this tree.
  Writing $\Length{\EMST{A}}$ for the length of the EMST of $A$, etc., we therefore have
  \begin{align}
    \norm{\Codomain{0}{f}}_1 &= \tfrac{1}{2} \Length{\EMST{A}} ; 
      \label{eqn:MST} \\
    \norm{\Domain{0}{f}}_1 &= \tfrac{1}{2} \Length{\EMST{A_0}} + \tfrac{1}{2} \Length{\EMST{A_1}} + \omega_0 ,
      \label{eqn:MST01}
  \end{align}
  in which $\omega_0 < \infty$ is the cutoff used for the death of the one $0$-cycle that does not die at any finite radius.\footnote{We consider reduced homology, nevertheless there can be points at infinity caused by considering the disjoint union of two sets, $A_0\sqcup A_1$.
  In particular, such points appear in $\Domain{0}{}$, $\Kernel{0}{}$, and $\Relative{1}{}$.}
  Equation~\eqref{eqn:MST} implies $\Ccodomain{0} = \frac{1}{2} \Cmst$, and since $A_0$ and $A_1$ are two stationary Poisson point processes with half the intensity each, Equation~\eqref{eqn:MST01} implies $\Cdomain{0} = {\sqrt{2} c}/{4} + {\sqrt{2} c}/{4} = {\sqrt{2} c}/{2}$, where we ignore the contribution of the cutoff, which in any case is a constant and thus of lower order than $\sqrt{n}$.
  Since $A = A_0 \cup A_1$, every $0$-cycle in the codomain is the image of a $0$-cycle in the domain.
  In addition, they are all born at radius $0$ so the image diagram is equal to the codomain diagram.
  Together with relation \eqref{eqn:short1} for $p = 0$, this implies $\Cimage{0} = \frac{1}{2} \Cmst$ and $\Ckernel{0} = \frac{1}{2} ( \sqrt{2} - 1 ) \Cmst$.

  \smallskip \noindent
  \textsc{Degree-$1$ domain and codomain diagrams.}
  In the alpha complex, a $1$-cycle is given birth by a critical edge and death by a critical triangle.
  The expected sums of radii are $\sqrt{n}$ for the critical edges and $\frac{3}{4} \sqrt{n}$ for the critical triangles \eqref{eqn:L21trianglestoo}.
  All critical triangles give death, but only some critical edges give birth, while the others give death. The sum of the radii of the latter is half the length of the EMST.
  This implies $\Expect{\norm{\Codomain{1}{f}}_1} = \Expect{F_2} - \Expect{F_1} + \tfrac{1}{2} \Expect{\Length{\EMST{A}}}$,
  and therefore $\Ccodomain{1} = \frac{1}{2} \Cmst - \frac{1}{4}$.
  Furthermore, $\Cdomain{1} = \sqrt{2} (\frac{1}{2} \Cmst - \frac{1}{4} )$, since $A_0$ and $A_1$ are stationary Poisson point processes of half the intensity each.

  \smallskip \noindent
  \textsc{Remaining diagrams.}
  By definition of $\CLmst$, we have $\Crelative{1} = \frac{1}{2} \CLmst$.
  There are now sufficiently many constants expressed in terms of $\Cmst$ and $\CLmst$ so that the remaining constants follow using the linear relations for the $1$-norms implied by the short exact sequences that connect the corresponding homology groups.
  Using Equations~\eqref{eqn:short3}, \eqref{eqn:short2}, \eqref{eqn:short1}, in this sequence, we get
  \begin{align}
    \Ccokernel{1} &= \Crelative{1} - \Ckernel{0}
                 = \tfrac{1}{2} \CLmst - \tfrac{1}{2} (\sqrt{2}-1) \Cmst ; \\
    \Cimage{1}    &= \Ccodomain{1} - \Ccokernel{1}
                 = \left( \tfrac{1}{2} \Cmst - \tfrac{1}{4} \right) - \left( \tfrac{1}{2} \CLmst - \tfrac{1}{2} (\sqrt{2}-1) \Cmst \right)
                 = \tfrac{1}{2} \sqrt{2} \Cmst - \tfrac{1}{4} - \tfrac{1}{2} \CLmst; \\
    \Ckernel{1}   &= \Cdomain{1} - \Cimage{1}
                 = \tfrac{1}{2} \sqrt{2} \left( \Cmst - \tfrac{1}{2} \right) - \tfrac{1}{2} \left( \sqrt{2} \Cmst - \tfrac{1}{2} - \CLmst \right)
                 = \tfrac{1}{2} \CLmst - \tfrac{1}{4} (\sqrt{2}-1) ,
  \end{align}
  and finally $\Crelative{2} = \Ckernel{1}$ because $\Ccokernel{2} = 0$.
\end{proof}

\subsection{Comparisons}
\label{sec:5.2}

\begin{figure}[hbt]
  \vspace{0.0in}
  \resizebox{!}{2.1in}{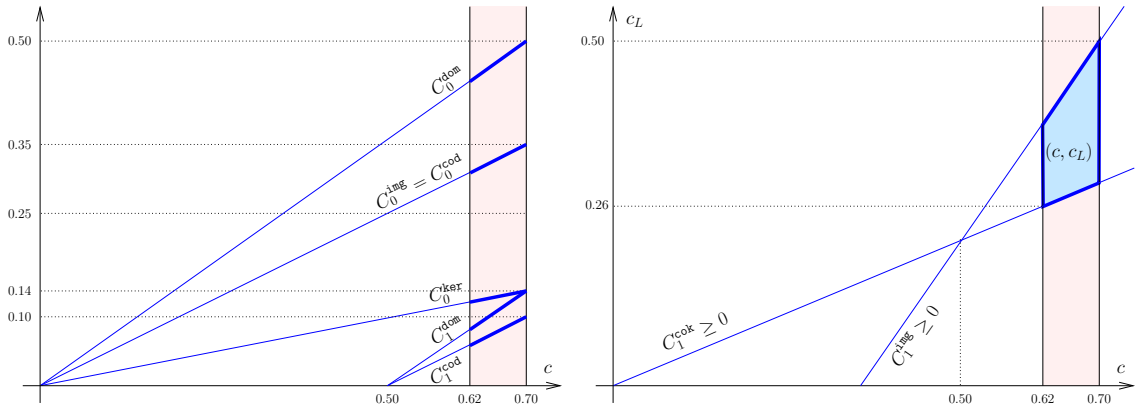}
    \vspace{-0.0in}
    \caption{\footnotesize \emph{Left panel:} the dependence of six constants on $\Cmst$, which lies somewhere between $0.6289$ and $0.7072$.
    Within this vertical strip, the affine maps are totally ordered.
    \emph{Right panel:}
    the zero-sets of the two constants that depend on $\Cmst$ as well as $\CLmst$.
    The non-negativity of these constants and the known bounds on $\Cmst$ imply the restriction of the possible pairs $(\Cmst, \CLmst)$ to the \emph{blue} quadrangular region.}
    \label{fig:Constants}
\end{figure}

All eleven non-zero constants in Table~\ref{tab:6-pack} depend either on $\Cmst$ or $\CLmst$ or both, so we do not know their exact values.
All constants are necessarily non-negative, since the $1$-norms of the diagrams are, and we get additional information from the known lower and upper bounds for $\Cmst$; see the left panel in Figure~\ref{fig:Constants}, which shows the affine maps of the six constants that depend on $\Cmst$ but not on $\CLmst$.
The non-negativity of $\Ccodomain{1}$ implies $\Cmst \geq \sfrac{1}{2}$, which is an alternative proof of Roberts' lower bound \cite{Rob68}.
Indeed, Roberts' argument is entirely different and uses a stochastic result about nearest neighbor pairs proved by Feller~\cite[page 14]{Fel66}.
Within the vertical strip between the known bounds for $\Cmst$, these maps are totally ordered, which implies
\begin{align}
  \Ccodomain{1} < \Cdomain{1} \leq \Ckernel{0} < \Cimage{0} = \Ccodomain{0} < \Cdomain{0} .
  \label{eqn:comparison1}
\end{align}
Consider the second inequality in \eqref{eqn:comparison1}, which requires Gilbert's upper bound of $\Cmst \leq \sfrac{\sqrt{2}}{2}$ \cite{Gil65}.
Since $\Cdomain{1}$ is equal to $\Ckernel{0}$ at this upper bound, and $\Cdomain{1}$ exceeds $\Ckernel{0}$ above the bound, an independent proof of this inequality would be an alternative proof of Gilbert's bound.
In this context, it is worth noting that these inequalities apply only to the expected case, and $\norm{\Domain{1}{}}_1 > \norm{\Kernel{0}{}}_1$ is quite possible in any particular case.
The same analysis for the three constants that depend on $\CLmst$ but not on $\Cmst$ leads to
\begin{align}
  \Ckernel{1} = \Crelative{2} < \Crelative{1} .
  \label{eqn:comparison2}
\end{align}
The non-negativity of $\Ckernel{1}$ and $\Crelative{2}$ implies $\CLmst \geq \frac{1}{2} (\sqrt{2}-1)$, but this bound can be improved.
\begin{theorem}
  \label{thm:lunar_asymptotic_constant}
  Let $\CLmst$ be the constant such that $\CLmst \sqrt{n}$ is the expected cost of the lunar EMST of $n$ points sampled uniformly at random in $[0,1]^2$ and randomly $2$-colored, in the limit, when $n$ goes to infinity.
  Then $0.2605 \leq \CLmst \leq 0.5$.
\end{theorem}
\begin{proof}
  Consider the two constants in Table~\ref{tab:6-pack} that depend on both asymptotic constants: $\Ckernel{1} = \frac{1}{2} \CLmst - \frac{1}{2} (\sqrt{2}-1) \Cmst$ and $\Cimage{1} = \frac{1}{2} \sqrt{2} \Cmst - \frac{1}{4} - \frac{1}{2} \CLmst$.
  Since the $1$-norms of the persistence diagrams are necessarily non-negative, we get
  \begin{align}
    (\sqrt{2}-1) \Cmst  \leq  \CLmst  \leq  \sqrt{2} \Cmst - \tfrac{1}{2} .
  \end{align}
  Plugging in the know bounds for $\Cmst$ gives the bounds for $\CLmst$.
  Indeed, plugging in $0.6289 \leq \Cmst$ on the left, we get $0.2605 \leq \CLmst$, and plugging in $\Cmst \leq \sfrac{\sqrt{2}}{2}$ on the right, we get $\CLmst \leq 0.5$.
\end{proof}

\section{Computational Experiments}
\label{sec:6}

This section presents the results of computational experiments aimed at estimating the expected $1$-norms of the persistence diagrams for the bi-chromatic setting explained in Section~\ref{sec:5}, and the related asymptotic constants $\Cmst$ and $\CLmst$.

\subsection{Experimental Set-up}
\label{sec:6.1}

We consider points selected uniformly at random in the unit square with and without periodic boundary conditions.
The difference between the two topologies becomes acute near the boundary of the square, where a point has possibly more near neighbors on the torus than in the square.
Hence, the minimum spanning tree on the torus is necessarily shorter than or at most as long as the minimum spanning tree for the same points in the square.
Nonetheless, the asymptotic constant, $\Cmst$, is the same for the two topologies, as proved in \cite{Jai93} and with different methods in \cite{Yuk98}.

\begin{figure}[hbt]
    \centering \vspace{-0.0in} 
  \subfloat[\footnotesize domain, $p = 0$.]{\label{fig:domain0}
    \centering
    \includegraphics[width=0.32\linewidth, keepaspectratio]{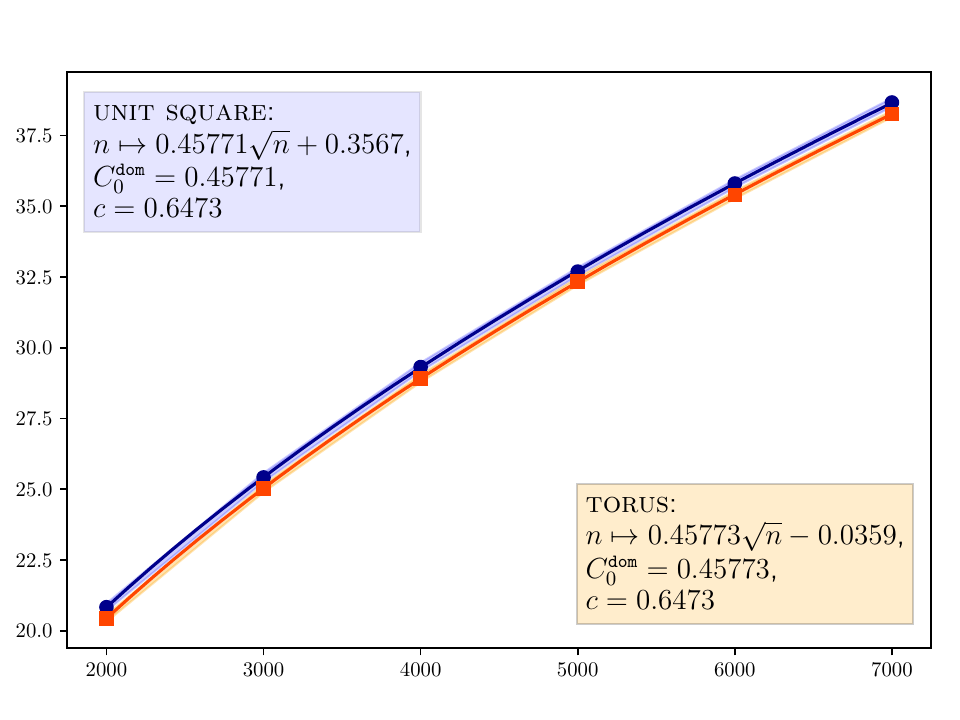}}
    \subfloat[\footnotesize image, $p = 0$.]{\label{fig:image0}
    \centering
    \includegraphics[width=0.32\linewidth, keepaspectratio]{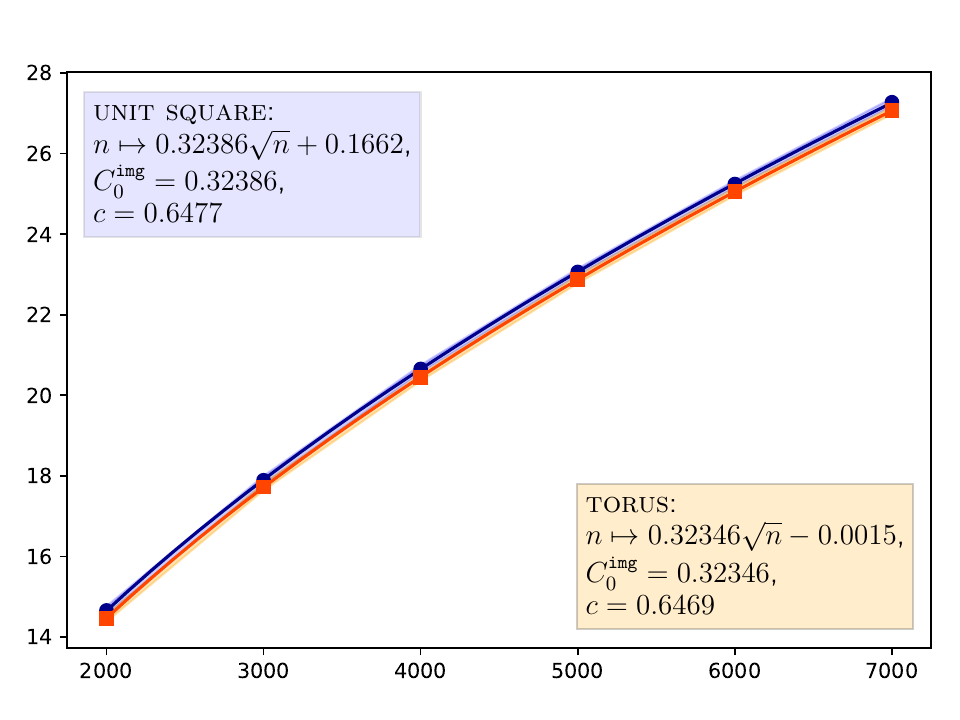}}
    \subfloat[\footnotesize codomain, $p = 0$.]{\label{fig:codomain0}
    \centering
    \includegraphics[width=0.32\linewidth, keepaspectratio]{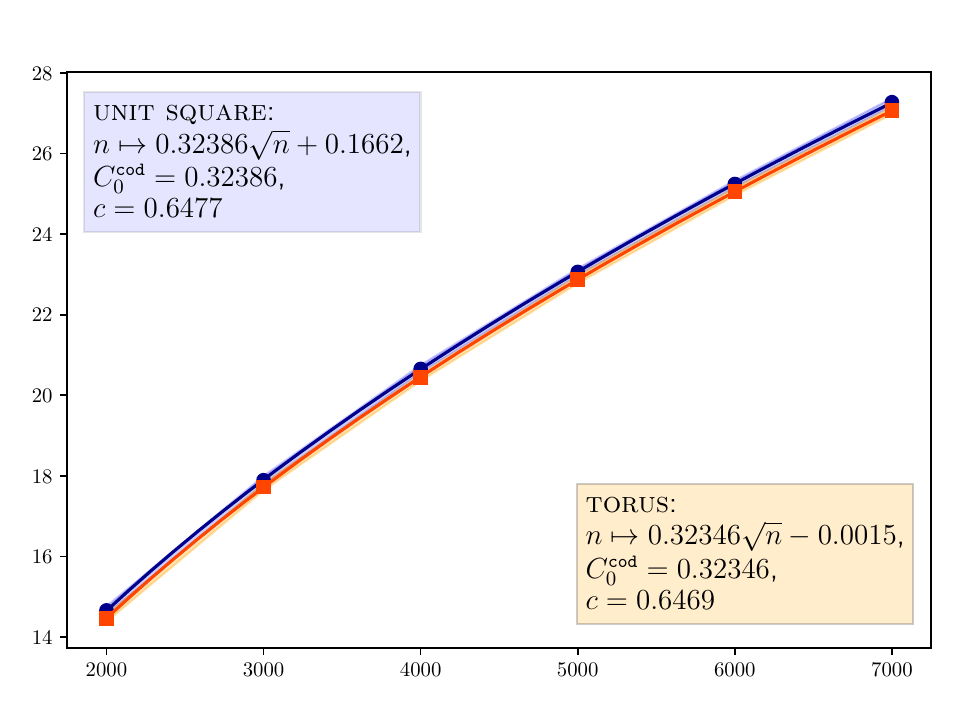}} \\
  \subfloat[\footnotesize kernel, $p = 0$.]{\label{fig:kernel0}
    \centering
    \includegraphics[width=0.32\linewidth, keepaspectratio]{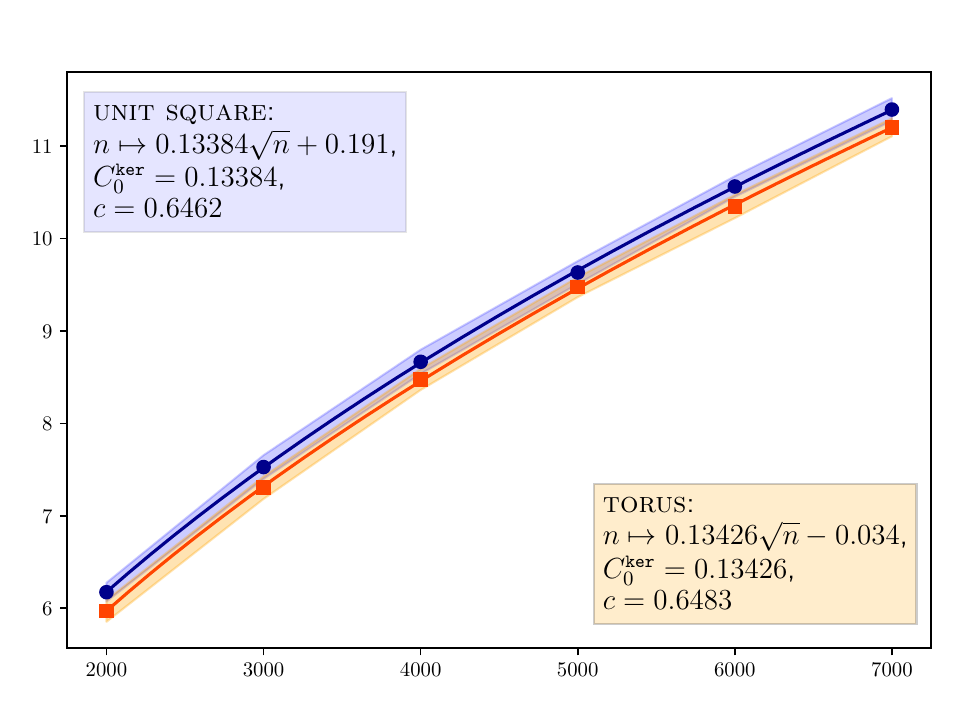}}
    \subfloat[\footnotesize relative, $p = 1$.]{\label{fig:relative1}
    \centering
    \includegraphics[width=0.32\linewidth, keepaspectratio]{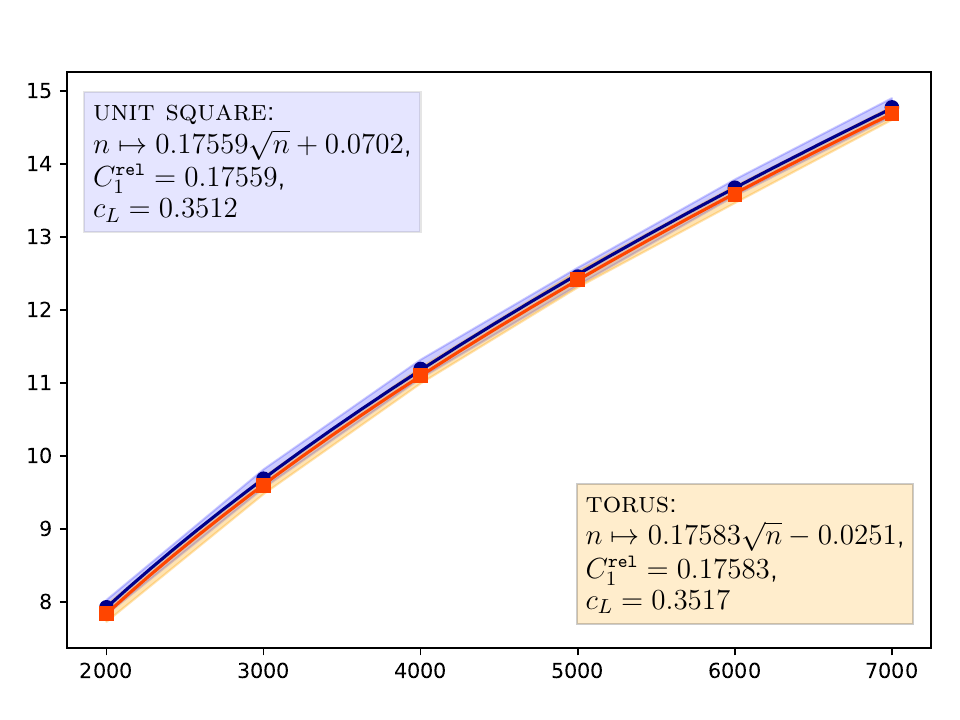}}
    \subfloat[\footnotesize cokernel, $p = 1$.]{\label{fig:cokernel1}
    \centering
    \includegraphics[width=0.32\linewidth, keepaspectratio]{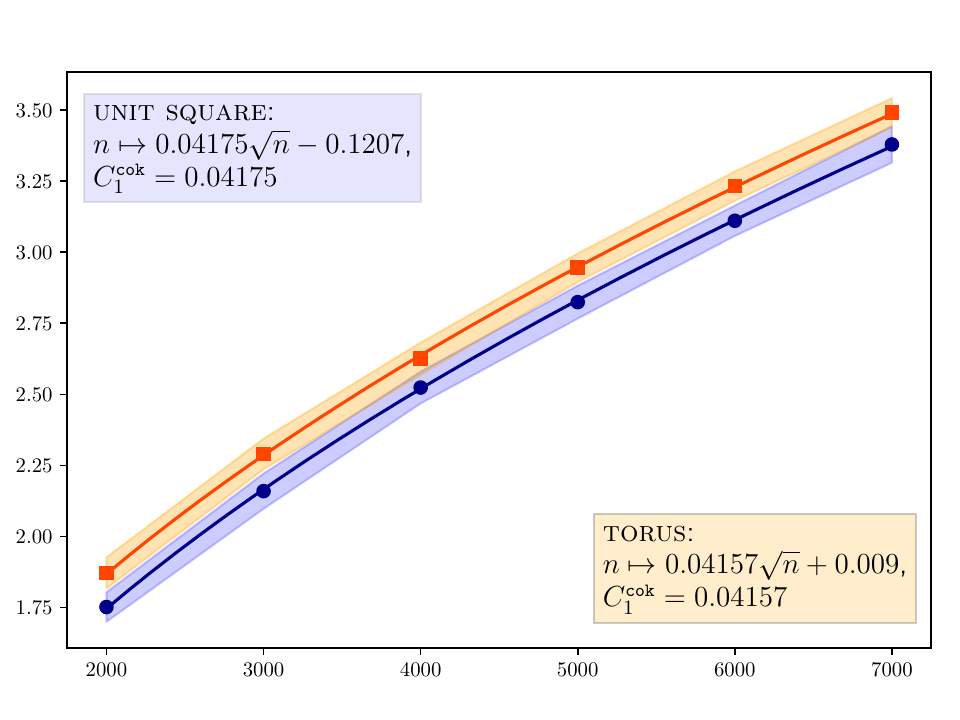}} \\
  \subfloat[\footnotesize domain, $p = 1$.]{\label{fig:domain1}
    \centering
    \includegraphics[width=0.32\linewidth, keepaspectratio]{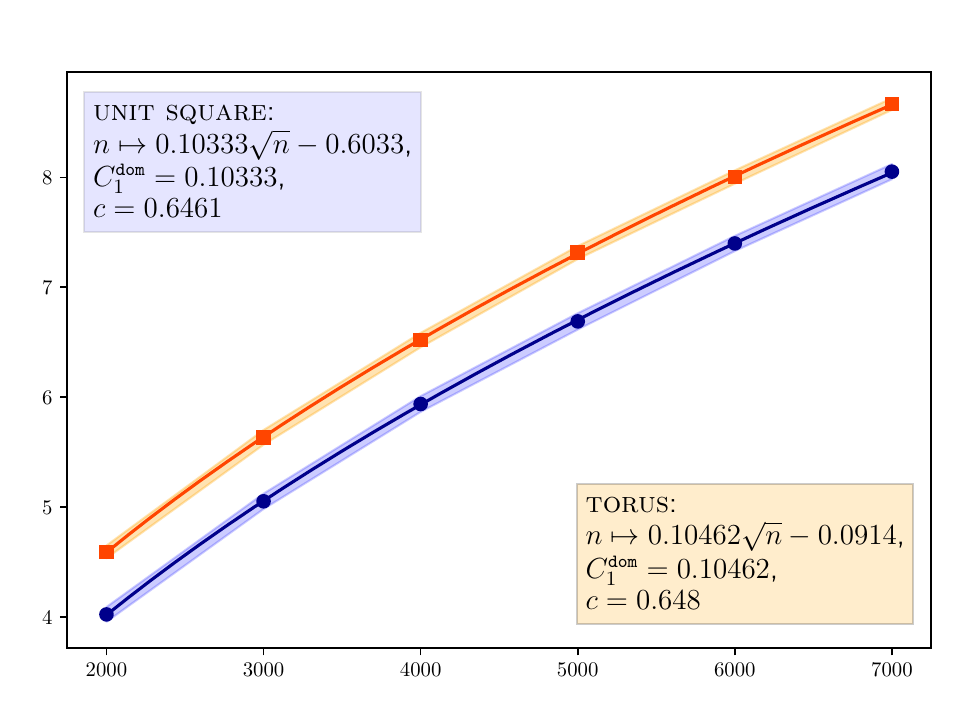}}
    \subfloat[\footnotesize image, $p = 1$.]{\label{fig:image1}
    \centering
    \includegraphics[width=0.32\linewidth, keepaspectratio]{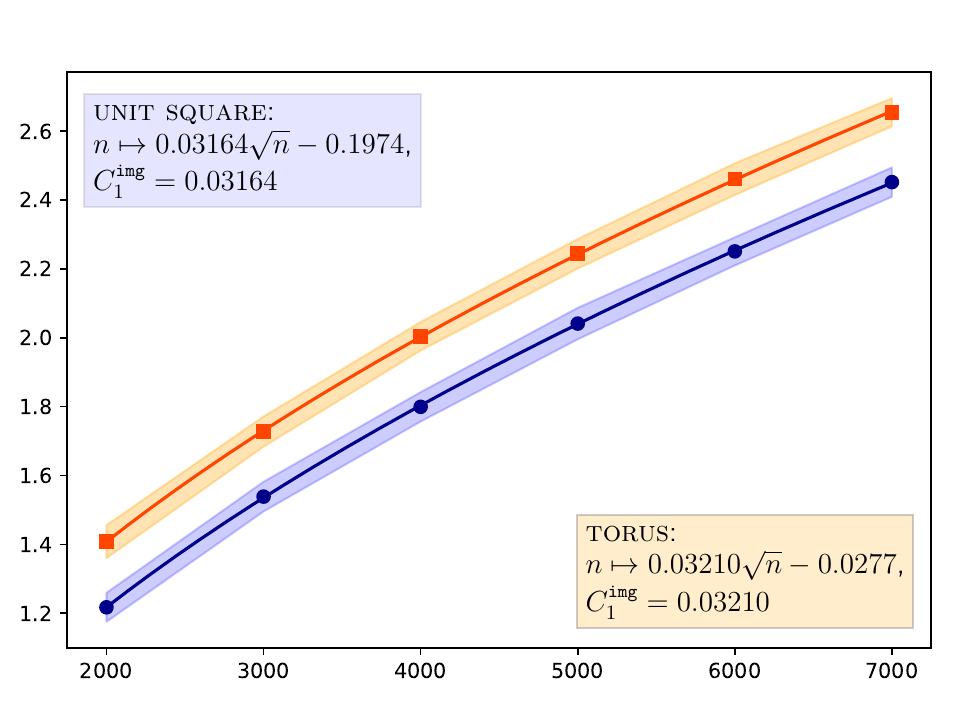}}
    \subfloat[\footnotesize codomain, $p = 1$.]{\label{fig:codomain1}
    \centering
    \includegraphics[width=0.32\linewidth, keepaspectratio]{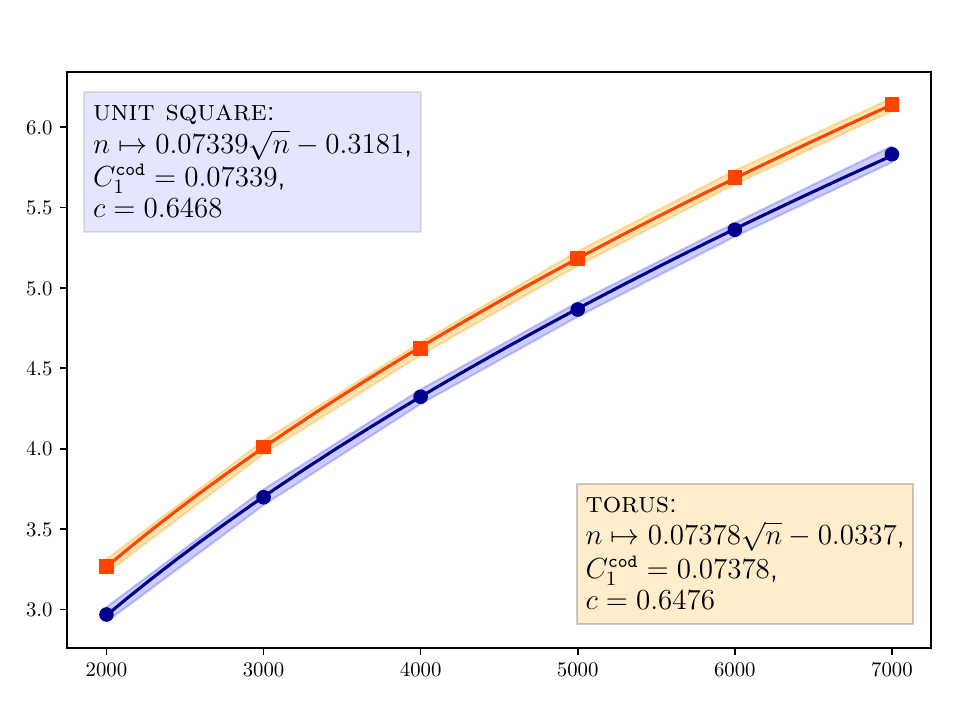}} \\
    \hspace{-2.10in}
    \subfloat[\footnotesize kernel, $p = 1$.]{\label{fig:kernel1}
    \centering
    \includegraphics[width=0.32\linewidth, keepaspectratio]{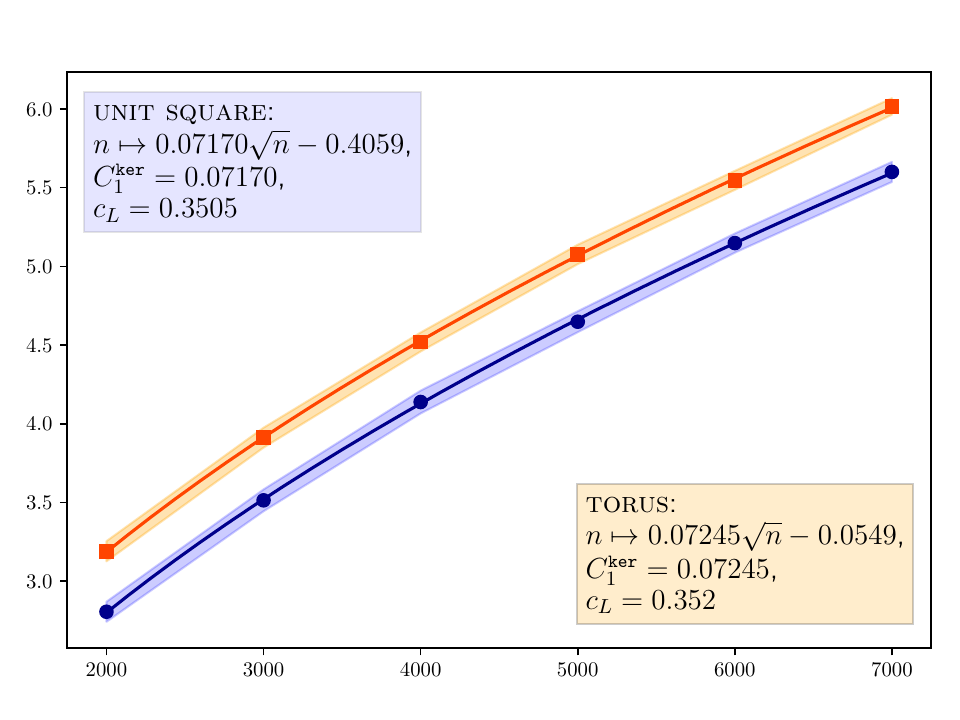}}
    \subfloat[\footnotesize relative, $p = 2$.]{\label{fig:relative2}
    \centering
    \includegraphics[width=0.32\linewidth, keepaspectratio]{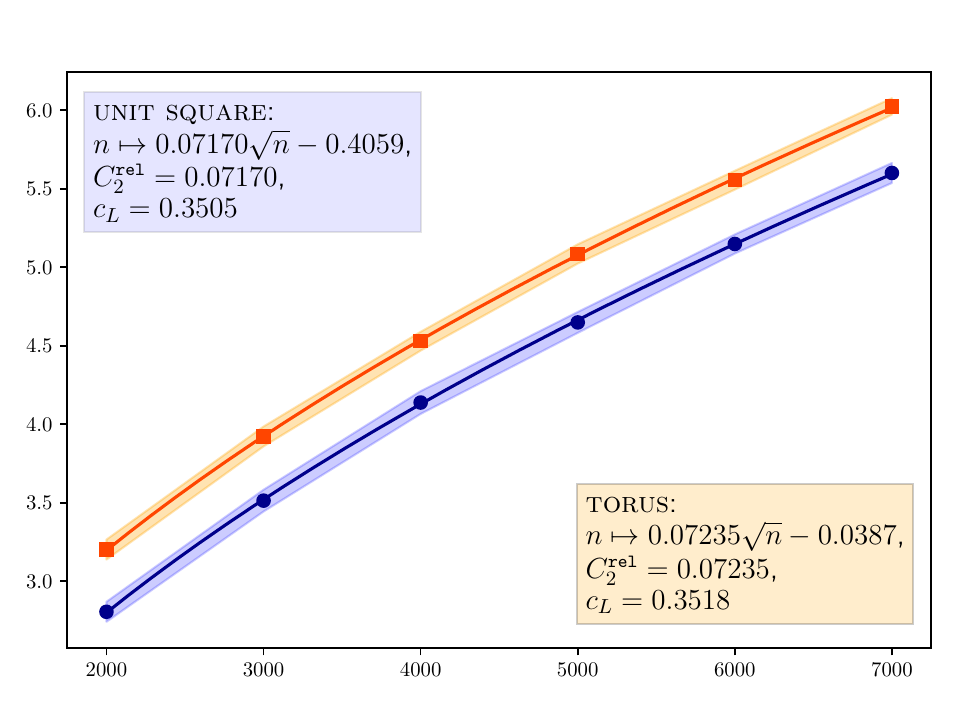}}
  \caption{\footnotesize The chromatic persistence diagrams organized in a snake-like fashion, \textsc{(c)-(b)-(a)-(d)-(e)-(f)-(i)-(h)-(g)-(j)-(k)}, with relations $\textsc{(c)=(b)}$, $\textsc{(a)=(b)+(d)}$, $\textsc{(e)=(d)+(f)}$, $\textsc{(i)=(f)+(h)}$, $\textsc{(g)=(h)+(j)}$, $\textsc{(k)=(j)}$ implied by short exact sequences.
  Each plot shows the best-fit curve for the average $1$-norms with the number of points increasing from \emph{left to right} (\emph{blue} for unit square and \emph{orange} for torus). 
  As illustrated by the colored strips along the fitting curves, the standard deviation is consistently small.
  In six of the eleven cases, we calculate the implied estimate of $\Cmst$, and in three the implied estimate of $\CLmst$.
  }
  \label{fig:plots}
\end{figure}

\smallskip
In a single experiment, we sample $n$ points in the unit square, color each point $0$ or $1$ with equal probability, and compute the $1$-norms of the chromatic persistence diagrams for both topologies.
There are eleven diagrams of interest: the domain and codomain diagrams for degrees $p = 0, 1$, the relative diagram of the pair for $p = 1, 2$, the kernel and image diagrams for $p = 0, 1$, and the cokernel diagram for $p = 1$.
All other diagrams are either empty or contain only information about the essential cycles (the homology of the square or the torus), which we ignore.
For each of these eleven cases, we compute the average $1$-norm of one hundred experiments each, for $n$ between $2000$ and $7000$, find the best-fit curve of the form $a_1 \sqrt{n} + a_0$, and interpret $a_1$ as an estimate of the corresponding asymptotic constant. 

\subsection{Findings}
\label{sec:6.2}
As predicted, the estimated length of the minimum spanning tree is slightly less for the torus than the square, but only by a very small margin.
As a consequence, the $1$-norms of the domain, image, and codomain diagrams in degree-$0$ are slightly smaller for the torus than the square.
We notice a somewhat larger difference for the sum of radii of the critical triangles, which is responsible for the more substantial deviation of the estimated $1$-norms of the degree-$1$ diagrams as well as the degree-$2$ relative diagram.
Curiously, the estimate of $\Ccokernel{1}$ for the square is slightly larger than that for the torus, while the best-fit curves show the opposite order; see panel \textsc{(f)} in Figure~\ref{fig:plots} and similar reversals in \textsc{(a)}, \textsc{(d)}, and \textsc{(e)}.
Short of a better explanation, we believe that the difference is too small to be a reliable signal.

\smallskip
Our intuition that the torus topology is a better approximation of the (infinite) Poisson point process in $\Rspace^2$ does not have strong support from the data.
For example, we may compare the estimates of $\Crelative{1} - \Crelative{2}$ with the theoretically predicted value of $(\sqrt{2}-1) \tfrac{1}{4} \approx 0.10355$; see Theorem~\ref{thm:constants}.
These estimates are $0.17559 - 0.07170 = 0.10389$ for the square and $0.17583 - 0.07235 = 0.10348$ for the torus, which straddle the predicted value with a small advantage for the latter.
For further comparison, we have repeated the experiments for a Poisson point process in the unit square and got virtually the same results (data not shown).

\smallskip
In summary, we venture that the asymptotic constant of the EMST length is between $0.646$ and $0.648$, and that of the lunar EMST cost is between $0.350$ and $0.352$, which is between the lower and upper bounds mentioned in Sections~\ref{sec:3.1} and \ref{sec:5.2}.

\section{Discussion}
\label{sec:7}

The main result of this paper is the improvement of the best lower bound for the asymptotic constant of the EMST length in two dimensions from $0.6008$ to $0.6289$.
We also conduct computational experiments, which suggest that the constant is between $0.646$ and $0.648$; compare with the current best upper bound of $0.7072$.
Our interest in the asymptotic constant is motivated by its relation to the $1$-norms of the chromatic persistence diagrams that arise in the study of randomly $2$-colored points in the Euclidean plane.
Purely algebraic considerations give rise to linear relations between the $1$-norms of these diagrams \cite{CDES23}, and we use these relations to cast light on the asymptotic constant from different directions.
A particular feature of these relations raises the following question:
\begin{itemize}
  \item Can the current best upper bound of Gilbert~\cite{Gil65} be improved below $\sfrac{\sqrt{2}}{2}$?
  It is curious that at this particular value we have $\Ckernel{0} = \Cdomain{1}$, as seen in the left panel of Figure~\ref{fig:Constants}.
  A proof that the expected $1$-norm of the degree-$1$ domain diagram cannot exceed that of the degree-$0$ kernel diagram would therefore give an alternative argument for Gilbert's upper bound.
\end{itemize}
There are several obvious directions for extending the results of this paper:  to maps different from the mono-chromatic to the bi-chromatic complexes, to three and more colors, to norms of the persistence diagrams different from the $1$-norm, and to three and higher dimensions.

\subsection*{Acknowledgements}
The authors thank Anton Nikitenko for his constructive criticism of an earlier version of this paper.

\newpage

\end{document}

%% file: Figs-length/Constants.pdf_t
\begin{picture}(0,0)%
\includegraphics{Figs-length/Constants.pdf}%
\end{picture}%
\setlength{\unitlength}{3947sp}%
\begingroup\makeatletter\ifx\SetFigFont\undefined%
\gdef\SetFigFont#1#2#3#4#5{%
  \reset@font\fontsize{#1}{#2pt}%
  \fontfamily{#3}\fontseries{#4}\fontshape{#5}%
  \selectfont}%
\fi\endgroup%
\begin{picture}(19229,7004)(-311,452)
\put(-149,6764){\makebox(0,0)[rb]{\smash{{\SetFigFont{14}{16.8}{\rmdefault}{\mddefault}{\updefault}{\color[rgb]{0,0,0}$0.50$}%
}}}}
\put(-149,4964){\makebox(0,0)[rb]{\smash{{\SetFigFont{14}{16.8}{\rmdefault}{\mddefault}{\updefault}{\color[rgb]{0,0,0}$0.35$}%
}}}}
\put(-149,1964){\makebox(0,0)[rb]{\smash{{\SetFigFont{14}{16.8}{\rmdefault}{\mddefault}{\updefault}{\color[rgb]{0,0,0}$0.10$}%
}}}}
\put(-149,2414){\makebox(0,0)[rb]{\smash{{\SetFigFont{14}{16.8}{\rmdefault}{\mddefault}{\updefault}{\color[rgb]{0,0,0}$0.14$}%
}}}}
\put(8401,539){\makebox(0,0)[b]{\smash{{\SetFigFont{14}{16.8}{\rmdefault}{\mddefault}{\updefault}{\color[rgb]{0,0,0}$0.70$}%
}}}}
\put(6001,539){\makebox(0,0)[b]{\smash{{\SetFigFont{14}{16.8}{\rmdefault}{\mddefault}{\updefault}{\color[rgb]{0,0,0}$0.50$}%
}}}}
\put(7426,539){\makebox(0,0)[b]{\smash{{\SetFigFont{14}{16.8}{\rmdefault}{\mddefault}{\updefault}{\color[rgb]{0,0,0}$0.62$}%
}}}}
\put(7351,1139){\rotatebox{27.0}{\makebox(0,0)[rb]{\smash{{\SetFigFont{20}{24.0}{\rmdefault}{\mddefault}{\updefault}{\color[rgb]{0,0,0}$\Ccodomain{1}$}%
}}}}}
\put(7351,1889){\rotatebox{33.0}{\makebox(0,0)[rb]{\smash{{\SetFigFont{20}{24.0}{\rmdefault}{\mddefault}{\updefault}{\color[rgb]{0,0,0}$\Cdomain{1}$}%
}}}}}
\put(7351,2414){\rotatebox{11.0}{\makebox(0,0)[rb]{\smash{{\SetFigFont{20}{24.0}{\rmdefault}{\mddefault}{\updefault}{\color[rgb]{0,0,0}$\Ckernel{0}$}%
}}}}}
\put(7351,4664){\rotatebox{27.0}{\makebox(0,0)[rb]{\smash{{\SetFigFont{20}{24.0}{\rmdefault}{\mddefault}{\updefault}{\color[rgb]{0,0,0}$\Cimage{0} = \Ccodomain{0}$}%
}}}}}
\put(7351,6239){\rotatebox{35.0}{\makebox(0,0)[rb]{\smash{{\SetFigFont{20}{24.0}{\rmdefault}{\mddefault}{\updefault}{\color[rgb]{0,0,0}$\Cdomain{0}$}%
}}}}}
\put(18301,539){\makebox(0,0)[b]{\smash{{\SetFigFont{14}{16.8}{\rmdefault}{\mddefault}{\updefault}{\color[rgb]{0,0,0}$0.70$}%
}}}}
\put(15901,539){\makebox(0,0)[b]{\smash{{\SetFigFont{14}{16.8}{\rmdefault}{\mddefault}{\updefault}{\color[rgb]{0,0,0}$0.50$}%
}}}}
\put(17326,539){\makebox(0,0)[b]{\smash{{\SetFigFont{14}{16.8}{\rmdefault}{\mddefault}{\updefault}{\color[rgb]{0,0,0}$0.62$}%
}}}}
\put(9751,6764){\makebox(0,0)[rb]{\smash{{\SetFigFont{14}{16.8}{\rmdefault}{\mddefault}{\updefault}{\color[rgb]{0,0,0}$0.50$}%
}}}}
\put(9826,3914){\makebox(0,0)[rb]{\smash{{\SetFigFont{14}{16.8}{\rmdefault}{\mddefault}{\updefault}{\color[rgb]{0,0,0}$0.26$}%
}}}}
\put(17776,4814){\makebox(0,0)[b]{\smash{{\SetFigFont{20}{24.0}{\rmdefault}{\mddefault}{\updefault}{\color[rgb]{0,0,0}$(\Cmst, \CLmst)$}%
}}}}
\put(10801,1439){\rotatebox{24.0}{\makebox(0,0)[lb]{\smash{{\SetFigFont{20}{24.0}{\rmdefault}{\mddefault}{\updefault}{\color[rgb]{0,0,0}$\Ccokernel{1} \geq 0$}%
}}}}}
\put(14851,1139){\rotatebox{55.0}{\makebox(0,0)[lb]{\smash{{\SetFigFont{20}{24.0}{\rmdefault}{\mddefault}{\updefault}{\color[rgb]{0,0,0}$\Cimage{1} \geq 0$}%
}}}}}
\put(8776,1064){\makebox(0,0)[b]{\smash{{\SetFigFont{20}{24.0}{\rmdefault}{\mddefault}{\updefault}{\color[rgb]{0,0,0}$\Cmst$}%
}}}}
\put(18676,1064){\makebox(0,0)[b]{\smash{{\SetFigFont{20}{24.0}{\rmdefault}{\mddefault}{\updefault}{\color[rgb]{0,0,0}$\Cmst$}%
}}}}
\put(10126,7139){\makebox(0,0)[lb]{\smash{{\SetFigFont{20}{24.0}{\rmdefault}{\mddefault}{\updefault}{\color[rgb]{0,0,0}$\CLmst$}%
}}}}
\put(-149,3764){\makebox(0,0)[rb]{\smash{{\SetFigFont{14}{16.8}{\rmdefault}{\mddefault}{\updefault}{\color[rgb]{0,0,0}$0.25$}%
}}}}
\end{picture}%